% Logic Eprints
%Submitted 1008 Fri May 20, 1994 by: andreas.blass@math.lsa.umich.edu (andreas blass )
%logic/blass/scc
%

%\input amstex
\magnification=1200
\documentstyle{amsppt}

\define\bp{\boldsymbol\Pi^0_}
\define\bs{\boldsymbol\Sigma^0_}
\define\bP{\boldsymbol\Pi^1_}
\define\bS{\boldsymbol\Sigma^1_}
\define\lp{ \Pi^0_}
\define\ls{ \Sigma^0_}
\define\lP{ \Pi^1_}
\define\lS{ \Sigma^1_}
\redefine\aa{{\bold a}}
\define\bb{{\bold b}}
\define\cc{{\bold c}}
\define\ii{{\bold i}}
\define\rr{{\bold r}}
\redefine\ss{{\bold s}}
\define\ttt{{\bold t}}
\define\dd{{\bold d}}
\define\Meager{\text{Meager}}
\define\NWD{\text{NWD}}
\define\Null{\text{Null}}

\topmatter
\title Simple Cardinal Characteristics \\
of the Continuum
\endtitle
\rightheadtext{Simple Cardinal Characteristics}
\author Andreas Blass \endauthor
\address
Mathematics Dept., University of Michigan, Ann Arbor, MI 48109, U.S.A.
\endaddress
\email
ablass\@umich.edu
\endemail
\thanks
Partially supported by NSF grant DMS 88-01988.
\endthanks
 \abstract We classify many cardinal characteristics of the continuum
according to the complexity, in the sense of descriptive set theory,
of their definitions.  The simplest characteristics
($\boldsymbol\Sigma^0_2$ and, under suitable restrictions,
$\boldsymbol\Pi^0_2$) are shown to have pleasant properties, related
to Baire category.  We construct models of set theory where
(unrestricted) $\boldsymbol\Pi^0_2$-characteristics behave quite chaotically
and no new characteristics appear at higher complexity levels.  We
also discuss some characteristics associated with partition theorems
and we present, in an appendix, a simplified proof of Shelah's theorem
that the dominating number is less than or equal to the independence
number.
\endabstract
\endtopmatter
\document

\head
1. Introduction
\endhead

      Cardinal characteristics of the continuum are cardinal numbers, 
 usually between  $\aleph_1$   and  $\cc = 2^{\aleph_0}$    inclusive, 
that give information about 
 the real line  $\Bbb R$  or the closely related sets  
$\Cal P(\omega)$ (the power set of the set  $\omega$  of natural numbers),
$[\omega]^\omega$  (the set of infinite subsets of $\omega$), 
and ${}^\omega\omega$ (the set of functions from  $\omega$ to $\omega$).  We give a few examples here 
 (others will be given later) and refer to [3, 23] for more examples and an 
 extensive discussion.

      The most obvious characteristic is  $\cc$, 
the cardinality of  $\Bbb R$  and of 
 the other sets,  $\Cal P(\omega)$, etc., mentioned above.
      
Baire category gives rise to several characteristics, of which we 
 mention here the covering number 
$$
      {\bold {cov}}(B) = 
       \text{ minimum number of meager sets needed to cover } \Bbb R 
 $$
and the  uniformity number
$$
      {\bold {unif}}(B) = 
       \text{ minimum cardinality of a non-meager set of reals.}
 $$
Lebesgue measure gives rise to analogous characteristics, 
${\bold {cov}}(L)$  and  ${\bold {unif}}(L)$,
  defined by putting ``measure zero" in place of 
``meager" in the previous definitions.

      For  $X,Y\in[\omega]^\omega$, we say that  $X$  {\sl splits\/}  $Y$  
if both  $X\cap Y$  and  $Y - X$  are infinite.  
Then the {\sl splitting number\/}  $\ss$  is the minimum 
 cardinality of a splitting family, 
i.e., a family  $\Cal S\subseteq[\omega]^\omega$  such that 
 every  $Y\in[\omega]^\omega$   is split by some  $X\in\Cal S$.  
The {\sl unsplitting number\/} $\rr$ (sometimes called the reaping 
number or the refining number) is the minimum cardinality of an 
 unsplittable family, i.e., a family  $\Cal R\subseteq [\omega]^\omega$
 such that no single  $X\in[\omega]^\omega$
 splits all the sets  $Y\in\Cal R$.

      For  $f,g\in{}^\omega\omega$, we say that  $g$  
{\sl eventually majorizes\/}  $f$, written 
$ f <^*  g $, if  $f(n) < g(n)$  for all but finitely many  $n$.  
The {\sl dominating 
 number\/} $\dd$  is the minimum cardinality of a dominating family, i.e., a 
 family  $\Cal D\subseteq{}^\omega\omega$  such that 
every  $f\in{}^\omega\omega$  is eventually majorized by some 
 $g\in\Cal D$.  The {\sl bounding number\/}  $\bb$  
is the minimum cardinality of an unbounded 
 family, i.e., a family  $\Cal B\subseteq{}^\omega\omega$
  such that no single  $g\in{}^\omega\omega$  eventually 
 majorizes all the functions  $f\in\Cal B$.

      For  $X,Y\in[\omega]^\omega$, we say that  $X$  is 
{\sl almost included\/} in  $Y$, written 
 $X\subseteq^*Y$, if  $X - Y$  is finite.  
The tower number  $\ttt$  is the smallest ordinal 
(necessarily a regular cardinal) such that there is a $\ttt$-sequence 
 $(X_\alpha)_{\alpha<\ttt}$ from $[\omega]^\omega$ that is 
almost decreasing (i.e., $X_\alpha\subseteq^*X_\beta$ for $\beta<\alpha$) 
 and cannot be extended, i.e., no $Y\in[\omega]^\omega$
  satisfies  $Y\subseteq ^*X_\alpha$  for all $\alpha < \ttt$.
 Such a sequence, not necessarily of minimum length, 
is called a {\sl tower\/}.

      Many cardinal characteristics have definitions of the form 
``the minimum cardinality of a set  $\Cal X$  
of reals such that, for every real $y$, there 
 is at least one  $x\in\Cal X$ such that $R(x,y)$'', where  $R$ is some binary 
 relation on reals.  If, as is customary in set theory, we apply the name 
 reals not only to members of  $\Bbb R$ 
 but also to members of  $\Cal P(\omega)$ and   ${}^\omega\omega$, 
 then the definitions of  $\ss$, $\rr$, $\dd$, and  $\bb$
  are clearly of this form.  The 
 definitions of the covering and uniformity numbers 
for category and measure 
 can be put into this form by coding meager or measure zero Borel sets by 
 reals.  The definition of  $\ttt$, however, seems not to be of this form, since 
 the family  $\Cal X = \{X_\alpha\mid \alpha<\ttt\}$  
is subject to the additional requirement of being 
 well-ordered by  $\supseteq^*$.  
Other characteristics, particularly  ${\bold h}$  and  ${\bold g}$  (see 
 [23]) are even further from the simple form we are discussing.

      The primary purpose of this paper is to discuss the possibility of 
 classifying those characteristics definable in this simple form by 
 applying, to the relations  $R$  involved in these definitions, the familiar 
 hierarchical classifications of descriptive set theory.  Specifically, we 
 define for any pointclass  $\Gamma$
  (e.g. for any level of the Borel hierarchy) 
 a uniform $\Gamma$-characteristic to be 
an uncountable cardinal definable in the 
 simple form described above, with  $R \in \Gamma$.  
We also define a somewhat 
 broader notion of (non-uniform) $\Gamma$-characteristic,
 which encompasses  $\ttt$  and 
 several other familiar characteristics.  Our goal is to infer combinatorial 
 information about a characteristic  $\kappa$  
from the descriptive set-theoretic 
 information that  $\kappa$  is a $\Gamma$-characteristic 
or a uniform $\Gamma$-characteristic 
 for a reasonably small  $\Gamma$.

      For example, when  $\Gamma = \bs2$, 
we can show that Martin's axiom implies 
 that all  $\bs2$-characteristics equal  $\cc$.  
In fact, in  ZFC  alone, all 
$\bs2$-characteristics are  $\geq{\bold {cov}}(B)$, 
which equals  $\cc$  under  MA  and which is 
 itself a uniform  $\bp1$-characteristic.

      For  $\bp2$, the situation is not so pleasant.  
$\aleph_1$   is a uniform 
 $\bp2$-char\-ac\-ter\-is\-tic, and there is a great deal of arbitrariness as to which 
 cardinals are  $\bp2$-characteristics and which are not.  Nevertheless, for 
 certain well-behaved  $\bp2$-characteristics, 
we can obtain results dual to those for  $\bs2$.

      This paper is organized as follows.  Section 2 contains definitions, 
 examples, and elementary properties of $\Gamma$-characteristics 
and uniform $ \Gamma$-characteristics.  
Sections 3 through 5 contain the results cited above, 
 $\bs2$  being treated in Section 3 and  $\bp2$   in Sections 4 and 5.  
In Section 6, 
 we record some information about characteristics related to Ramsey's 
 theorem and other partition theorems.  These characteristics first 
 attracted my attention as interesting (and possibly new) examples of  
$\bs2$ 
 and  $\bp2$   characteristics,
but they have since acquired considerable combinatorial interest 
in their own right.  Finally, following the example of 
 Vaughan [23], we present in an appendix a proof of Shelah's theorem that 
 the independence number  $\ii$  is no smaller than  $\dd$.  This proof is a 
 reformulation of Shelah's original proof [23] avoiding a few unnecessary 
 complications.  Although it contains no new (vis \`a vis [23]) ideas, its 
 relative simplicity seems to justify recording it in print.

\head
2. Definitions, Basic Properties, and Examples
\endhead

      Throughout this paper,  $\Gamma$  denotes a pointclass 
in the sense of  descriptive set theory.  
The specific  $\Gamma$'s  that we consider will usually 
 be low levels of the Borel hierarchy, 
particularly  $\bp1$   and  $\bp2$, 
their lightface analogs $\lp1$   and  $\lp2$, 
and the class  $OD\Bbb R$  of relations 
 ordinal-definable from reals.  
The relevance of  $OD\Bbb R$  is  that it includes
 all ``reasonable'' pointclasses, so, by proving that a cardinal is not an 
 $OD\Bbb R$-characteristic, we establish that it is not a 
$\Gamma$-characteristic for 
 ``any''  $\Gamma$.
 
\demo{Definition}
An uncountable cardinal $\kappa$ is a 
$\Gamma$-{\sl characteristic\/} if there is 
 a family of $\kappa$ sets, each in  $\Gamma$, 
such that   ${}^\omega\omega$ is covered by the family 
 but not by any subfamily of cardinality  $<\kappa$.  
An uncountable cardinal $\kappa$
 is a {\sl uniform $\Gamma$-characteristic\/} 
if there is a binary relation  $R$  on  ${}^\omega\omega$
 such that  $R\in\Gamma$  and such that $\kappa$ is the minimum cardinality of a family 
 $\Cal X\subseteq{}^\omega\omega$  
such that for every  $y\in{}^\omega\omega$
 there exists  $x\in\Cal X$  with  $R(x,y)$.
\enddemo

        The remainder of this section is devoted to examples and elementary 
 properties of $\Gamma$-characteristics and uniform 
$\Gamma$-characteristics.  
The following proposition is obvious.
\proclaim{Proposition 1}
 (a)\quad  If  $\Gamma\subseteq\Delta$ then every 
(uniform) $\Gamma$-characteristic is a 
 (uniform) $\Delta$-characteristic.\newline
      (b)  Every  $\Gamma$-characteristic  lies between  
$\aleph_1$   and  $\cc$, inclusive.\newline
      (c)  If  $\Gamma$  is closed under pre-images 
by continuous functions, then 
 every uniform  $\Gamma$-characteristic is a $\Gamma$-characteristic. 
\qed 
\endproclaim 

      In connection with (c), we note that the only continuous functions 
 needed are  $y\mapsto(x,y)$  for arbitrary fixed  $x$.  
The hypothesis of (c) 
 means in practice that  $\Gamma$  is a boldface pointclass; 
for lightface pointclasses  $\Gamma$  
there are no $\Gamma$-characteristics but there are uniform 
 $\Gamma$-characteristics.
      
In defining (uniform) $\Gamma$-characteristics, we used the descriptive 
 set-theorist's usual version  ${}^\omega\omega$ of the ``reals''.  
Had we used the actual 
 reals $\Bbb R$  or  $\Cal P(\omega)$  or  $[\omega]^\omega$
   instead, any (uniform) $\Gamma$-characteristic in 
 the modified sense would also be a (uniform) $\Gamma$-characteristic in the 
 original sense, provided  $\Gamma$  is closed under pre-images by continuous 
 functions with recursive codes.  This follows immediately from the 
well-known fact [18, p. 12] that ${}^\omega\omega$ can be mapped 
onto each of  $\Bbb R$, $\Cal P(\omega)$, and $[\omega]^\omega$  
by continuous functions with recursive codes.  
The specific $\Gamma$'s that we deal with will have this closure property, 
so, when verifying that some  cardinal is a 
(uniform) $\Gamma$-characteristic, we may use  
$\Bbb R$  or  $\Cal P(\omega)$  or  $[\omega]^\omega$
or similar spaces instead of  ${}^\omega\omega$.

     The next proposition lists some examples of 
$\Gamma$-characteristics, mostly uniform ones. 
 The proof in each case consists of exhibiting the family of 
 sets or the binary relation required by the definition and then verifying 
 membership in  $\Gamma$.  We shall see later that $\bs2$ can always be improved to $\bp1$,
 so part (b) of the proposition is not optimal.  In fact, it is easy to 
 prove directly that  $\dd$, $\rr$, and  ${\bold {cov}}(B)$  
are uniform  $\bp1$-characteristics. 
 But for the time being, we give the complexity bounds that follow directly 
 from the definitions.
 \proclaim{Proposition 2}
(a)  $\cc$  is a uniform  $\lp1$-characteristic.\newline
 (b)  $\dd$, $\rr$, ${\bold {cov}}(B)$, and  ${\bold {unif}}(L)$
  are uniform $\ls2$-characteristics.\newline
 (c)  $\bb$, $\ss$, ${\bold {unif}}(B)$, ${\bold {cov}}(L)$, and  
   $\aleph_1$
   are uniform  $\lp2$-characteristics.\newline
 (d)  $\ttt$  is a  $\bp2$-characteristic.
\endproclaim
 
 \demo{Proof}
(a)  Take  $R(x,y)$  to be  $x = y$, a $\lp1$ relation.

 (b)  For  $\dd$, take
$$
     R(x,y) \iff y \leq^* x \iff \exists k\,(\forall n \geq k)\, y(n) \leq x(n). 
$$
 For  $\rr$, work in $[\omega]^\omega$ and take
$$\align
   R(x,y) &\iff y  \text{ does not split } x \\
          &\iff \exists k\,[(\forall n \geq k)\,(n\in x \implies n\in y)\\  
  &\qquad\qquad\text{ or }  (\forall n\geq k)\,(n\in x \implies n\notin y)]. 
\endalign$$
For  ${\bold {cov}}(B)$, we need a coding of the meager sets, 
or rather of the 
 countable unions of nowhere dense closed sets, by ``reals''.  
Any nowhere dense closed subset $F$  of ${}^\omega\omega$
 is of the form
$$
       \{x\in {}^\omega\omega\mid
       \text{ No finite initial segment of } x \text{ is in } D\}
$$ 
 where  $D\subseteq{}^{<\omega}\omega $
(the set of finite sequences from  $\omega$) is dense in the sense 
 that every  $s\in {}^{<\omega}\omega$ has an extension
  $s{}^\frown t \in D$ (where  ${}^\frown$ means concatenation).  
If  $f:{}^{<\omega}\omega\to{}^{<\omega}\omega$,
 then  $\{s{}^\frown f(s)\mid s\in{}^{<\omega}\omega\}$  is such a  $D$, 
and every dense  $D$  includes one of this form.  Thus,
$$\align
      &\NWD(f) =\\
&\qquad \{x\in{}^\omega\omega\mid
(\forall s\in{}^{<\omega}\omega)\, s^\frown f(s)  
\text{ is not an initial segment of }x\}
\endalign $$
is a nowhere dense closed set, and every nowhere dense set is included in 
 one of this form.  Therefore, if  
$g:\omega\times{}^{<\omega}\omega\to{}^{<\omega}\omega$, then
$$\align
  &\Meager(g) = \\
&\qquad\{x\in{}^\omega\omega\mid
\exists n\,(\forall s\in{}^{<\omega}\omega) s^\frown g(n,s) 
\text{ is not an initial segment of }x \}
\endalign$$
 is a meager  $F_\sigma$   set, 
and every meager set is included in one of this form. 
 It follows that  ${\bold {cov}}(B)$  is the minimum size of a family
 $\Cal X$  of  $g$'s  (in 
${}^{\omega\times{}^{<\omega}\omega}({}^{<\omega}\omega)$
which is recursively homeomorphic to ${}^\omega\omega$)  
such that for every 
 $y\in{}^\omega\omega$ there is  $g\in\Cal X$  with  
$y\in \Meager(g)$.  By inspection of the definition 
of $\Meager(g)$, we see that the relation  $y\in\Meager(g)$  is a  
$\ls2$  
 relation  $R(g,y)$.  Therefore  ${\bold {cov}}(B)$  is a uniform  
$\ls2$-characteristic. 

      (Remark:  By regarding  ${\bold {cov}}(B)$  as the 
minimum number of nowhere dense 
 (rather than meager) sets needed to cover the reals, we could work with 
 $\NWD(f)$  rather than $\Meager(g)$  and get, with less work, the better result 
 that  ${\bold {cov}}(B)$  is a uniform $\lp1$-characteristic.  
But we shall need  $\Meager(g) $
 later, and the improvement from $\ls2$ to $\lp1$ will be automatic 
when we establish Proposition 3(a) below.)

      For  ${\bold {unif}}(L)$, we also need a coding,
 this time of the measure zero  $G_\delta$
 sets, by reals.  The idea is to code a measure zero  
$G_\delta$   set  $N$  by coding 
 a sequence of open sets  $U_n$, with measures  $\mu(U_n)\leq2^{-n}$, and 
with intersection  $N$; an open set  $U_n$   in turn is coded 
by listing codes for 
 basic open sets whose union is  $U_n$.  
For convenience, we work in   ${}^\omega2 $
(or equivalently  $\Cal P(\omega)$), 
where every finite sequence  $s\in{}^{<\omega}2$  determines a 
 basic open set
$$
  B_s=\{x\in{}^\omega2\mid s\text{ is an initial segment of }x\} 
 $$
of measure  $2^{-\text{length}(s)}$.  
We intend to use functions  
$g:\omega\times\omega\to{}^{<\omega}2$  to 
 code  $G_\delta$ sets  
$\bigcap_{n\in\omega}\bigcup_{k\in\omega}B_{g(n,k)}$, 
but we wish to ensure that only  $G_\delta$ 
 sets of measure zero are obtained.  Therefore, we set
$$\align
\Null(g) =& \{x\in{}^\omega\omega\mid
   [\forall n\,\exists k\,g(n,k)\text{ is an initial segment of }x]
\text{ and}\\
   &[\forall n\,\forall k\,\sum_{j=0}^k2^{-\text{length}(g(n,j))}\leq2^{-n}]\}.
\endalign $$
The second clause here means that  
$\mu(\bigcup_{j\in\omega}B_{g(n,j)})\leq2^{-n}$   
and so the coded $G_\delta$ set must have measure zero.  
If this clause, which is independent of  $x$,
 is not satisfied by  $g$, then  $\Null(g)=\emptyset$.  
Thus, $\Null(g)$  is always a 
$G_\delta$ set of measure zero, 
and every set of measure zero is included in one of 
 this form.  Therefore,  ${\bold {unif}}(L)$  is 
the minimum cardinality of a family  $\Cal X$ of reals 
(in ${}^\omega2$) such that, 
for every  $g$ (in ${}^{\omega\times\omega}({}^{<\omega}2)$ 
which is homeomorphic
to ${}^\omega\omega$),  there is  $x\in\Cal X$  with $x\notin\Null(g)$.  
The relation $x\notin\Null(g)$ is, by inspection of 
the definition of $\Null(g)$, $\ls2$, and so 
 ${\bold {unif}}(L)$  is a uniform  $\ls2$-characteristic.
 
(c)  For  $\bb$, take,
$$
R(x,y) \iff x\not\leq^*y\iff\forall k\,(\exists n\geq k)\,y(n)<x(n). 
 $$
For  $\ss$, take
$$\align
    R(x,y) &\iff x\text{ splits }y\\
              &\iff \forall k\,
                 [(\exists n\geq k)\,(n\in x \land n\in y) \text{ and }\\
                  &\qquad\qquad(\exists n\geq k)\,(n\notin x\land n\in y)].
\endalign$$
For ${\bold {unif}}(B)$, take
$$
                         R(x,g) \iff x \notin\Meager(g). 
$$
For ${\bold {cov}}(L)$, take
$$
                          R(g,y) \iff y\in\Null(g). 
$$
Finally, for $\aleph_1$, take, for  $x\in{}^\omega\omega$  and  
$y\in{}^{\omega\times\omega}\omega$,
$$\align
R(x,y)&\iff x\text{ differs from all the functions }
                  y(n,-)\text{ for }n\in\omega\\
          &\iff (\forall n)\,(\exists k)\,x(k)\neq y(n,k).
\endalign$$
All these $R$'s are $\lp2$ by inspection.
 
(d)  Let  $(X_\alpha)_{\alpha<\ttt}$ be as in the definition of  $\ttt$. 
For each  $\alpha<\ttt$, let
$$
  Q_\alpha=\{Y\in[\omega]^\omega\mid Y-X_\alpha\text{ is infinite}\}. 
 $$
Then $Q_\alpha$ is $\bp2$ for each $\alpha$, 
and the family $\{Q_\alpha\mid\alpha<\ttt\}$  covers $[\omega]^\omega$
 because the sequence $(X_\alpha)_{\alpha<\ttt}$
 cannot be extended.  A subfamily of 
 cardinality $<\ttt$  has the form $\{Q_\alpha\mid\alpha\in S\}$
  where  $S$  is a subset of $\ttt$  of 
 cardinality  $<\ttt$  and is therefore not cofinal in  $\ttt$. 
So let $\beta<\ttt$  be larger than every $\alpha\in S$.  
Then, for $\alpha\in S$, we have $X_\beta\subseteq^*X_\alpha$.  Therefore  
$X_\beta\notin\bigcup_{\alpha\in S}Q_\alpha$, 
and the subfamily fails to cover $[\omega]^\omega$.
\qed\enddemo

      We remark that the proof of (d) shows that, if a regular cardinal $\kappa$
 is the length of a tower, then $\kappa$ is a $\bp2$-characteristic.

\proclaim{Proposition 3}
(a)  Every uniform $\ls{n+1}$-char\-ac\-ter\-is\-tic is a uniform 
 $\lp n$-char\-ac\-ter\-is\-tic.\newline
 (b)  Every (uniform) $\bs{n+1}$-characteristic is a (uniform) 
$\bp n$-char\-ac\-ter\-is\-tic.\newline
(c)  Every $\bS2$-characteristic is a ${\boldsymbol\Delta}^1_1$-characteristic.
\endproclaim

 \demo{Proof}
(a) Let  $R$  be a $\ls{n+1}$ binary relation on  ${}^\omega\omega$, 
given by
$$
                         R(x,y) \iff (\exists n\in\omega)\,Q(n,x,y) 
 $$
where  $Q$  is $\lp n$.  If $\Cal X$  is a subset of ${}^\omega\omega$ 
such that $\forall y\,(\exists x\in\Cal X)\,R(x,y)$, 
 then $\omega\times\Cal X$  is a subset $\Cal Z$  of  
$\omega\times{}^\omega\omega$  such that  
$\forall y\,(\exists(n,x)\in\Cal Z)\,Q(n,x,y)$. 
 Conversely, for any such 
$\Cal Z\subseteq\omega\times{}^\omega\omega$, 
we obtain such an $\Cal X\subseteq{}^\omega\omega$  by 
 taking $\Cal X=\{x\mid(\exists n)\,(n,x)\in\Cal Z\}$.  
In each case $\Cal X$ and $\Cal Z$ have the same 
 cardinality provided they are infinite.  So the uniform 
 $\ls{n+1}$-characteristic defined by  $R$  equals the uniform 
 $\lp n$-characteristic 
 defined by $Q\subseteq(\omega\times{}^\omega\omega)
\times{}^\omega\omega$. 
 
(b)  The uniform case is proved exactly like (a); just make all the $\Sigma$'s 
 and $\Pi$'s boldface.  The non-uniform case is (at least) equally easy; in the 
 given cover of ${}^\omega\omega$ by $\bs{n+1}$ sets, 
replace each of these sets  $A$  with 
 countably many $\bp n$ sets whose union is  $A$.
 
(c)  Every $\bS2$ set is a union of $\aleph_1$ Borel 
(i.e., $\boldsymbol\Delta^1_1$) sets  [18,p.96]. 
 So if $\kappa$ is a $\bS2$-characterstic and 
$\Cal X$  is a covering of ${}^\omega\omega$ by $\kappa$ $\bS2$
 sets such that no subfamily of size $<\kappa$ 
covers ${}^\omega\omega$, we can replace each 
 of the $\bS2$ sets in $\Cal X$  by its 
$\aleph_1$ Borel constituents to obtain a new 
 covering $\Cal X'$  of ${}^\omega\omega$ by $\kappa$
(since $\kappa\geq\aleph_1$) $\boldsymbol\Delta^1_1$ sets such that 
no subfamily of size  $<\kappa$ covers  ${}^\omega\omega$. 
\qed\enddemo

      Neither Proposition 3 nor the list of examples in Proposition 2 
 mentioned $\bs1$-characteristics, and for a good reason.  There are none.  If
 ${}^\omega\omega$ is covered by a family of $\bs1$ 
(i.e., open) subsets, then it is covered  by a countable subfamily 
(i.e., the space ${}^\omega\omega$  is Lindel\"of), but we 
 required characteristics to be uncountable.
 
\proclaim{Corollary 4}
 $\dd$, $\rr$, ${\bold {cov}}(B)$, and ${\bold {unif}}(L)$  are 
uniform $\lp1$-char\-ac\-ter\-is\-tics.
\endproclaim

 \demo{Proof}
Combine Propositions 2(b) and 3(a). 
\qed\enddemo
  
    In view of Proposition 3, the $\Sigma$ levels of the arithmetical and 
 finite Borel hierarchies are, in the context of characteristics, equivalent 
 to the immediately preceding $\Pi$ levels.  (``Finite'' is inessential here. 
 Proposition 3(b) remains true, with the same proof, if  $n$  is allowed to be 
 transfinite.)  For the rest of the paper, we shall use the $\Pi$ rather than 
 the $\Sigma$ class.

      To avoid leaving an obvious and unnecessary gap in our list of 
 examples, we comment on the additivity and cofinality characteristics of 
 measure and category.  These characteristics are defined by
$$\align
{\bold{add}}(B) =&\text{minimum number of meager sets 
              whose union is not meager},\\
{\bold{cof}}(B) =&\text{minimum number of meager sets in a family}\\
              &\text{such that every meager set has a superset in the family},
\endalign$$
 and analogous definitons for  ${\bold{add}}(L) $ and  ${\bold{cof}}(L)$. 
In these definitions, 
 we can restrict attention to meager $F_\sigma$ sets and 
measure zero  $G_\delta$   sets, 
 i.e., to sets easily coded by reals.  The coding shows that these are 
 characteristics in the sense of our definitions and in fact that both 
 additivities are uniform $\lS1$-characteristics and both cofinalities are 
 uniform $\lP1$-characteristics.  
One can do considerably better by using combinatorial 
 characterizations of these cardinals.  Specifically, by writing out the 
 characterizations of  ${\bold{add}}(L)$  and  ${\bold{cof}}(L)$
 implicit in Theorems 0.9 and 
 0.10 of [1] (and using Lemma 0.5 to eliminate a quantifier), we find that 
 ${\bold{add}}(L)$ is a uniform $\lp2$-characteristic and  
${\bold{cof}}(L)$  is a uniform 
$\ls2$- (hence $\lp1$-)characteristic.  
(Actually, Theorems 0.9 and 0.10 are 
concerned with the equations  ${\bold{add}}(L) = \cc$ 
and  ${\bold{cof}}(L) = \cc$, but the 
generalizations we need can be proved exactly the same way.  The same remark applies to other results cited below.)  
For the Baire category characteristics, we 
 have that  ${\bold{add}}(B)$  is either  ${\bold {cov}}(B)$  or 
$\bb$, whichever is smaller [15], 
 and  ${\bold{cof}}(B)$  is either ${\bold {unif}}(B)$  or 
$\dd$, whichever is larger [16].  So both 
 ${\bold{add}}(B)$  and  ${\bold{cof}}(B)$ are uniform $\lp2$-characteristics.
      We remark that Bartoszy\'nski's elegant characterization [1, Thm. 1.7] 
 of  ${\bold {cov}}(B)$  as the smallest cardinality of any 
$\Cal X\subseteq{}^\omega\omega$  such that 
 $$
(\forall y\in{}^\omega\omega)\,(\exists x\in\Cal X)\,
(\exists n)\,(\forall k\geq n) x(k)\neq y(k)
$$ 
gives an alternative proof that 
${\bold {cov}}(B)$  is a uniform $\ls2$- 
(hence $\lp1$@-)char\-ac\-ter\-is\-tic.

\head
3.  Lower Bound for $\bp1$-Characteristics
\endhead

\proclaim{Theorem 5}
 If $\kappa$ is a $\bp1$-characteristic, then $\kappa\geq{\bold {cov}}(B)$. 
\endproclaim  

\demo{Proof}
Let  $\Cal X$  be a covering of ${}^\omega\omega$  by $\kappa$
 closed (i.e., $\bp1$) sets.   We 
 shall show that either $\kappa\geq{\bold {cov}}(B)$  or a countable
subfamily of  $\Cal X$  covers  ${}^\omega\omega$.   
In particular, if  $\Cal X$ witnesses that $\kappa$
 is a $\bp1$-characteristic, then 
 the second alternative is impossible, so the theorem will follow. 

Fix a countable base for the topology of ${}^\omega\omega$, 
and let  $U$  be the 
 union of all the basic open sets  $B$  such that some countable subfamily of 
 $\Cal X$ covers  $B$.   
(We do not claim that any such  $B$'s exist; $U$ might be 
 empty.)  As the base is countable, $U$ itself can be covered by a countable 
 subfamily of $\Cal X$.   So if  $U={}^\omega\omega$, we are done. 

Henceforth, we assume that $U\neq{}^\omega\omega$, 
and we let  $C={}^\omega\omega-U$.   
So  $C$  is a nonempty closed subset of ${}^\omega\omega$.   
We claim that  $C$  is perfect.  
 Indeed, if  $C$  had an isolated point  $x$, 
then some basic neighborhood  $B$  of  $x$  would be included in  
$U\cup\{x\}$, which can be covered by countably 
 many sets from  $\Cal X$ --- countably many to cover  $U$  
and one more to cover  $x$.  
 But then  $B\subseteq U$, contrary to  $x\in B\cap C = B - U$. 

So  $C$  is a perfect subset of  ${}^\omega\omega$, 
and clearly  $C$  is covered by the 
 closed sets  $X\cap C$  for  $X\in\Cal X$.   
We claim that each of these closed sets 
 is nowhere dense in $C$.   Indeed, if this claim were false, there would be 
 $X\in\Cal X$  and a basic open set  $B$  such that  
$B\cap C$  is nonempty and  $\subseteq X\cap C$.  
 But then  $B\subseteq U\cup X$, 
so  $B$  can be covered by a countable subfamily of 
 $\Cal X$ --- countably many to cover  $U$, plus  $X$.   
But then  $B\subseteq U$, contrary to 
 $B\cap C\neq\emptyset$. 

Thus, $C$  is covered by $\kappa$ nowhere dense (in $C$) sets $X\cap C$.   
But it is well-known that the Baire covering number 
for any perfect subset  $C$  of ${}^\omega\omega$ is the same 
as for ${}^\omega\omega$.   
So we must have $\kappa\geq{\bold {cov}}(B)$. 
\qed\enddemo

 \proclaim{Corollary 6}
All of  $\dd$, $\rr$, and ${\bold {unif}}(L)$  are 
$\geq{\bold {cov}}(B)$.\qed 
\endproclaim
      This corollary is, of course, well-known, but the usual proof of 
${\bold {unif}}(L)\allowmathbreak\geq{\bold {cov}}(B)$, 
due to Rothberger [20], uses more specific information 
 about Lebesgue measure.   Our proof, by contrast, 
uses only that  ${\bold {unif}}(L)$ 
 is a $\bs2$-\ (and therefore $\bp1$@-)characteristic.   
Specifically, if we replace 
 the ideal  $L$  of measure zero sets by any ideal that, like  $L$, has a 
 cofinal subfamily indexed by reals, say  
$\{S_x\mid x\in{}^\omega\omega\}$, such that the 
 relation $y\in S_x$ is $\bp2$ (or even $\bp2$  in  $x$  for each fixed  $y$) 
then ${\bold {unif}}$
 for this ideal is $\geq{\bold {cov}}(B)$.   
Rothberger's proof, on the other hand, is 
 symmetric between category and measure; it shows that  
${\bold {unif}}(B)\geq{\bold {cov}}(L) $
 for the same reason as  ${\bold {unif}}(L)\geq{\bold {cov}}(B)$.   
Our proof lacks this symmetry 
 because  ${\bold {unif}}(B)$  and  ${\bold {cov}}(L) $ are not 
$\bp1$-\ but $\bp2$-characteristics.   We 
 shall discuss this matter further in Section 5. 
 
\proclaim{Corollary 7}
Martin's axiom implies that every $\bp1$-characteristic 
is equal to $\cc$. 
\endproclaim

\demo{Proof}
Martin's axiom, even when weakened to apply only to countable (rather
 than ccc) partial orders, implies ${\bold {cov}}(B) = \cc$.  [5,14].
\qed\enddemo

\head
4. Consistency Results for $\bp2$ and Higher Characteristics
\endhead

      All the familiar cardinal characteristics of the continuum  --- those 
 defined in Section 1 above as well as numerous others [3,23] --- are equal 
 to  $\cc$  under Martin's axiom, and Corollary 7  may be viewed as a partial 
 explanation of this fact.  Unfortunately, this explanation does not extend 
 beyond the $\bp1$ level, since $\aleph_1$   is a uniform
$\lp2$-characteristic. 
 This leaves open several possibilities.  
For example, $\aleph_1$ might be the only 
 exception to a very general situation.  That is, there might be a large 
 pointclass  $\Gamma$  (perhaps even  $\Gamma$ = OD$\Bbb R$)  
such that, under  MA, all 
 $\Gamma$-characteristics greater than $\aleph_1$ are equal to  $\cc$.  
On the other hand, 
 one might think that, as one goes up the Borel hierarchy, more and more 
 small cardinals become characteristics; since $\aleph_1$ is a 
$\bp2$-characteristic, 
 perhaps $\aleph_2$ is a $\bp3$-characteristic and so on.

In this section, we present consistency results showing that none of 
 these possibilities is provable.  We shall see that, from $\bp2$ upward, 
 there is considerable arbitrariness in the characteristics.  For example, 
 as the following theorem shows, it is consistent that, for all $n\in\omega$,
$$\align
 \aleph_n\text{ is a }\bp2\text{-characteristic }
&\iff  \aleph_n\text{ is an OD$\Bbb R$-characteristic }\\
&\iff  n \text{ is a power of }17.
 \endalign$$

\proclaim{Theorem 8}
Assume  GCH, and let  $A$  be a subset of $\omega$ containing  1  but 
 not  0.  Then there is a forcing extension of the universe in which 
 $\cc=\aleph_{\omega+1}$; 
the $\bp2$-characteristics are $\aleph_n$ for $n\in A$,
$\aleph_\omega$, and $\aleph_{\omega+1}$; and 
 these are the only  OD$\Bbb R$-characteristics.
\endproclaim

      We shall prove a somewhat more general result, but it seems worthwhile 
 to point out first that the presence of $\aleph_\omega$ 
among the $\bp2$-char\-ac\-ter\-is\-tics 
 is unavoidable if  $A$  is infinite.  
Indeed, if  $\Gamma$  is any non-trivial pointclass closed 
 under pre-images by recursively coded continuous functions, 
then the supremum of 
 any countably many $\Gamma$-characteristics 
is also a $\Gamma$-characteristic.  Indeed, 
 if $\Cal X_n$ witnesses that $\kappa_n$ is a $\Gamma$-characteristic, 
then 
$\{\{(n)^\frown x\mid x\in X\}\mid n\in\omega, X\in\Cal X_n\}$
witnesses that  $\sup_n\kappa_n$ is a $\Gamma$-characteristic.

\proclaim{Theorem 9}
Assume  GCH, and let  $C$  be a closed set of uncountable 
 cardinals containing  $\aleph_1$, 
containing all uncountable cardinals $\leq|C|$ and 
 containing the immediate successors of all its members of 
cofinality  $\omega$. 
 Then there is a notion of forcing satisfying the countable chain condition 
 and forcing that $\cc=\max(C)$  and that both 
the set of $\bp2$-characteristics 
 and the set of OD$\Bbb R$-characteristics are equal to  $C$.
\endproclaim

      We remark that, because  $C$  is closed, it has a largest element, so 
 the equation $\cc=\max(C)$ makes sense.

Clearly, Theorem 8 is a special case of Theorem 9 with  
$C=\{\aleph_n\mid n\in A\}\cup\{\aleph_\omega,\aleph_{\omega+1}\}$.  
In contrast to the situation with Theorem 8, we do not know 
 that all the hypotheses in Theorem 9 are really needed in their full 
 strength.  We certainly need that $C$ be closed under $\omega$-limits 
(see the remark preceding the theorem), 
that it have a largest element, 
that this element not have cofinality $\omega$
(so that $\cc=\max(C)$ is consistent), and 
 that it contain $\aleph_1$ (see Proposition 2(a, c)).  
But the remaining hypotheses 
 might be mere artifacts of the proof technique.

 \demo{Proof}
The proof consists of first ensuring that every $\kappa\in C$  is a 
$\bp2$-characteristic by forcing a 
maximal almost disjoint family of $\kappa$ subsets of $\omega$  
and second showing that no cardinals $\lambda\notin C$  are 
 OD$\Bbb R$-characteristics in the forcing extension.  
The first part of the proof 
 uses Hechler's technique [9] for forcing maximal almost disjoint families 
 of different sizes.  The second is related to a theorem of Miller [17] 
 that, when many independent Cohen reals are added to a model of  GCH, no 
 cardinal strictly between $\aleph_1$ and $\cc$ is a Borel-characteristic; we 
 strengthen the conclusion from ``Borel'' to ``OD$\Bbb R$'', and we work with a 
 Hechler-type model rather than the Cohen model.

 \proclaim{Lemma 10}
If $\Cal X$ is a maximal almost disjoint family of infinite subsets 
 of  $\omega$, then  $|\Cal X|$ is a $\bp2$-characteristic.
\endproclaim

 \demo{Proof}
We work in $[\omega]^\omega$ and define, for each $X\in\Cal X$,
$$
 M(X) = \{Y\in[\omega]^\omega\mid X\cap Y\text{ is infinite}\}. 
 $$
The family $\{M(X)\mid X\in\Cal X\}$  covers $[\omega]^\omega$, 
by maximality of  $\Cal X$.  But no 
 proper subfamily $\{M(X)\mid X\in\Cal X'\}$ with  
$\Cal X'\subset\Cal X$  covers  $[\omega]^\omega$, because if
 $Y\in\Cal X-\Cal X'$  then  $Y\notin M(X)$  for any $X\in\Cal X'$.  
By counting quantifiers, we 
 see that each  $M(X)$  is a $\bp2$ set, 
so  $|\Cal X|$  is a $\bp2$-characteristic.  
\qed\enddemo

      Returning to the proof of the theorem, we recall that we wish to 
 force, for each $\kappa\in C$, 
a maximal almost disjoint family of cardinality $\kappa$; 
 by Lemma 10, this will ensure that every $\kappa\in C$  
is a $\bp2$-characteristic in 
 the extension (provided cardinals are preserved, which they will be).  
For each $\kappa\in C$, let $I_\kappa=\{(\kappa,\xi)\mid\xi<\kappa\}$, 
and let $I=\bigcup_{\kappa\in C}I_\kappa$. The maximal 
 almost disjoint family of size $\kappa$ that we adjoin 
will be indexed by $I_\kappa$, 
 so altogether we shall adjoin an $I$-indexed family of subsets of $\omega$. 
 Since we want the forcing to satisfy the countable chain condition, we use 
 finite conditions; we build the desired almost-disjointness into the 
 forcing, and genericity will yield the desired maximality.

      A forcing condition $p$ is a function into  2  whose domain is of the 
 form $F\times n$ where $F$  is a finite subset of $I$ 
and $n\in\omega$.  (We make the 
 customary identification of $n$  with  $\{0,1,\dots,n-1\}$.)  An extension of 
 $p:F\times n\to2$  is a condition  $p':F'\times n'\to2$ such that  
$p'\supseteq p$ 
(and therefore  $F'\supseteq F$  and $n'\geq n$) 
and, whenever $(\kappa,\xi)$ and $(\kappa,\eta)$  are 
 distinct elements of $I_\kappa$  (for the same $\kappa$)  
and $n\leq k < n'$, then 
 $p'(\kappa,\xi,k)$  and  $p'(\kappa,\eta,k)$ are not both  1.  
(Intuitively, we regard 
 $p:F\times n\to2$  as giving the following information about 
the generic sets $a_{\kappa,\xi}\subseteq\omega$ being adjoined.  
First, if  $p(\kappa,\xi,k)=1$ (resp. 0), then 
 $k\in a_{\kappa,\xi}$ (resp. $k\notin a_{\kappa,\xi}$), 
and second, if  $(\kappa,\xi)$ and $(\kappa,\eta)$  are distinct 
 elements of $F\cap I_\kappa$, then 
$a_{\kappa,\xi}\cap a_{\kappa,\eta}\subseteq n$.  Then the definition of 
 extension corresponds to giving more information.)  We call this notion of 
 forcing  $P$.

      Before proceeding, we should point out that this notion of forcing is 
 essentially a part of the forcing defined by Hechler in [9, Theorem 3.2]. 
 Hechler adjoins maximal almost disjoint families of all cardinalities from 
 $\aleph_1$ to $\cc$ and towers of all lengths 
$\leq\cc$  of uncountable cofinality, and 
 he works with a ground model where $\cc$  is (in the interesting cases) 
 already large so that his forcing does not alter cardinal exponentiation. 
 If one ignores the parts of Hechler's forcing that refer to towers and one 
 replaces the interval  $[\aleph_1,\cc]$  of cardinals 
(destined to become the sizes 
 of maximal almost disjoint families) by  $C$, then one obtains a notion of 
 forcing having a dense subset isomorphic to our  $P$.  The next two lemmas 
 are transcriptions for  $P$  of corresponding arguments in Hechler's proof; 
 we include their proofs for the reader's convenience.

 \proclaim{Lemma 11}
$P$  satisfies the countable chain condition.
\endproclaim

 \demo{Proof}
Let $\aleph_1$ elements $p_\alpha:F_\alpha\times n_\alpha\to2$  
of  $P$  be given ($\alpha<\aleph_1$).  By 
 passing to a subfamily of size $\aleph_1$, 
we can assume that all the  $n_\alpha$   are 
 the same  $n$  and that the $F_\alpha$ constitute a 
$\Delta$-system [11, p.225], i.e., 
 $F_\alpha\cap F_\beta$  is the same set  $K$  
for all  $\alpha<\beta<\aleph_1$.  Again, by passing to a 
 subfamily of size $\aleph_1$, we can assume that 
the restrictions $p_\alpha|K\times n$ 
 are all equal.  But then, for any  $\alpha<\beta$,
$p_\alpha\cup p_\beta$   is a common extension 
 of $p_\alpha$ and $p_\beta$.  
So the $\aleph_1$ given conditions do not form an antichain. 
\qed\enddemo

      (This proof shows that in any family of $\aleph_1$ conditions, there is a 
 subfamily of $\aleph_1$ conditions every finitely many of which 
have a common extension, i.e.,  $P$  has precaliber $\aleph_1$.)

      Let  $G$  be a $P$-generic filter over the universe  $V$.  (Formally, we 
 are passing to a Boolean-valued extension of  $V$.)  
For $(\kappa,\xi)\in I$  let
$$
a_{\kappa,\xi}=\{k\mid(\exists p\in G)\,p(\kappa,\xi,k)=1\}, 
 $$
and for $\kappa\in C$  let
$$
\Cal A_\kappa=\{a_{\kappa,\xi}\mid\xi<\kappa\}. 
$$

\proclaim{Lemma 12}
For each $\kappa\in C$, $\Cal A_\kappa$ is a maximal almost disjoint 
family of subsets of  $\omega$  in $V[G]$.
\endproclaim

\demo{Proof}
 Fix   $\kappa\in C$.  If  $\xi<\eta<\kappa$, then  $G$  contains 
some  $p: F\times n\to2$ 
 with both  $(\kappa,\xi)$  and  $(\kappa,\eta)$  in  $F$, 
because the set (in  $V$) of all such 
 $p$'s  is clearly dense.  Fix such a  $p\in G$  and consider an arbitrary 
 $q\in G$.  As  $G$  is a filter,  $p$  and  $q$ 
 have a common extension  $r$.  By 
 definition of extension, we cannot have  
$r(\kappa,\xi,k)=r(\kappa,\eta,k) = 1$  for any 
 $k\geq n$, and thus we cannot have  
$q(\kappa,\xi,k) = q(\kappa,\eta,k) = 1$ for any  $k\geq n$. 
 As  $q$  was arbitrary in  $G$, it follows that 
$a_{\kappa,\xi}\cap a_{\kappa,\eta}\subseteq n$.  So  $\Cal A_\kappa$   is 
 an almost disjoint family of subsets of  $\omega$.

 To prove maximality, suppose  $x$  were  (in $V[G]$) an infinite subset of 
 $\omega$  almost disjoint from $a_{\kappa,\xi}$     for all  $\xi<\kappa$.
  Because of the countable 
 chain condition (Lemma 11), $x$ has name  
$\dot x\in V$  that involves only countably 
 many conditions.  Fix a countable set  $J\subseteq I$ 
 such that all the conditions 
 involved in $\dot x$ have domains  $\subseteq J\times \omega$.  
Also fix a condition 
 $p:F\times n\to 2$  forcing ``$\dot x$  is an infinite 
subset of $\omega$ almost disjoint 
 from $\dot a_{\kappa,\xi}$     for all  $\xi < \kappa$'' 
(where $\dot a$  is the standard name for the 
 function  $a$  that sends  $(\kappa,\xi)$  to  $a_{\kappa,\xi}$, 
and where we have written  $\kappa$ 
 instead of its canonical name  $\check\kappa$.)  
Enlarging  $J$  if necessary, we assume 
 $F\subseteq J$.  Since  $\kappa\in C$, $\kappa$  
is uncountable, so fix  $\xi < \kappa$  with  $(\kappa,\xi)\notin J$. 
 Since  $p$  forces  ``$\dot x\cap\dot a_{\kappa,\xi}$ is finite'', 
it has an extension 
 $p': F' \times  n'\to 2$  forcing  
``$\dot x\cap\dot a_{\kappa,\xi}\subseteq m$''  
for some  $m\in\omega$.  Extending  $p'$ 
 further, we can assume  $n'\geq m$, and then, as  $m$  can be increased 
 trivially, we can assume  $m = n'$, so  $p'$  forces  
``$\dot x\cap\dot a_{\kappa,\xi}\subseteq n'$.''  Now, 
 as  $p'$  extends  $p$  and therefore forces  
``$\dot x$  is infinite and almost  disjoint from all $\dot a_{\kappa,\eta}$'', 
it also forces (since  $F'$  is finite) 
``there is $k\geq n'$  such that
$$
k\in\dot x \text{ and } (\forall\eta\in X)\, 
k\notin \dot a_{\kappa,\eta},\text{''} \tag1
$$
 where  $X = \{\eta\mid(\kappa,\eta)\in F'\cap J\}$.

      Extend  $p'$  to a condition  $q:  H \times  n''\to  2$  
forcing \thetag1 for a specific 
 $k\geq n'$.  Let  $F'' = H\cap J$  and let  $p'' = q\mid(F'' \times  n'')$.

 \proclaim{Claim}
   $ p''$  forces \thetag1 for the same value of  $k$.
\endproclaim

To prove this claim, suppose it failed, and extend  $p''$  to a condition 
 $r$  forcing the negation of \thetag1.  As  $p''$
  and all the conditions involved in 
 \thetag1 (i.e., in $\dot x$ and in $\dot a_{\kappa,\eta}$ 
for $\eta\in X$) have domains $\subseteq J\times\omega$, we can 
 take $r$  to also have domain $\subseteq J\times\omega$.  
Then the function  $q\cup r$  can be 
 extended by zeros to a condition that extends both  $q$  and  $r$.  This is 
 absurd, as  $q$  forces  \thetag1  and  $r$  forces its negation.  
So the claim is proved.

Notice that  $p''$  agrees with  $p'$  on the common part of their domains, 
 $(F'\cap J)\times n'$, because $q$  extends them both.  Extending  $p''$  if 
 necessary, we assume  $n'' > k$; then,
$$
(\forall\eta\in X)\, p''(\kappa,\eta,k) = 0 \tag2
$$
 because  $p''$  forces \thetag1.

      Define a function  $s:F'\times n''\to 2$  
by making  $s$  agree with $p'$  on 
 $F' \times  n'$, making  $s$  agree with  $p''$
  on  $(F'\cap J) \times  n''$, setting 
 $s(\kappa,\xi,k) = 1$, and setting all 
remaining values of  $s$  equal to zero.  It 
 is obvious that  $s$  extends $p'$ as a function; we claim it extends $p'$ as 
 a condition.  We must check that, for any  $(\lambda,\alpha)\neq (\lambda,\beta) $ both in  $F'\cap I_\lambda$ 
 and any  $j$  with  $n'\leq j<n''$, the values of $s$ at the two locations 
 $(\lambda,\alpha,j)$ and $(\lambda,\beta,j)$  are not both 1.  
If they were, then neither of these 
 values could be given by the first clause in the definition of $s$, because 
 $j\geq n'$ and the first clause gives values on  $F'\times n'$.  
Neither of these 
 values could be given by the last clause, since the last clause gives 
 values of zero.  One of the values could be given by the second clause, but 
 not both, for the values given by the second clause agree with the values 
 of  $q$, an extension of $p'$.  So one of the two values was given by the 
 third clause and the other by the second.  That is, $\lambda=\kappa$,
 $j = k$, and one 
 of $\alpha$ and $\beta$, say $\alpha$, is  $\xi$.  Then
$$
1 = s(\lambda,\beta,j) = s(\kappa,\beta,k) 
$$
 is given by the second clause, i.e., $(\kappa,\beta)\in F' \cap J$.  
But then
$$
s(\kappa,\beta,k) = p''(\kappa,\beta,k) = 0, 
$$
 by \thetag2.  This contradiction shows that $s$ extends $p'$ as a condition.

      Therefore,  $s$, like  $p'$, forces  ``$k\in\dot x$''  and 
 ``$\dot x\cap\dot a_{\kappa,\xi}\subseteq n'$.''  As 
 $k\geq n'$, $s$ must also force  ``$k\notin \dot a_{\kappa,\xi}$.''  
But this is absurd, as 
 $s(\kappa,\xi,k) = 1$.  This contradiction completes the proof of maximality.  
\qed\enddemo

      By the lemmas proved so far, every  $\kappa\in C$  is a  
$\bp2$-characteristic in 
 $V[G]$.  It remains to prove that in  $V[G]$ no cardinal  
$\lambda\notin C$  is an  OD$\Bbb R$-characteristic 
and that  $\cc = \max(C)$.  The latter actually follows from the 
 former, since  $\cc$  is the largest $\bp2$-characteristic, 
but it is also easy to 
 see directly, since the notion of forcing  $P$  has cardinality $\max(C)$, 
 whose cofinality is uncountable, and  $P$  satisfies the countable chain 
 condition, and  GCH  holds in the ground model  $V$.

      To complete the proof, consider any uncountable  
$\lambda\notin C$, and suppose 
 we have, in  $V[G]$, a $\lambda$-sequence of  OD$\Bbb R$  sets 
$X_\alpha  (\alpha <\lambda)$  that cover  ${}^\omega\omega$. 
 Fix a sequence of reals  $u_\alpha$  and a sequence of ordinals  
$\theta_\alpha$  such that  $X_\alpha$ 
 is ordinal-definable with real parameter  $u_\alpha$  
and in fact is the  $\theta_\alpha$th set 
 ordinal-definable from  $u_\alpha$ 
(in some standard well-ordering of the  OD($u_\alpha$)  sets).  
Choose in  $V$  names  $\dot X$, $\dot u$, $\dot\theta$ 
for the sequences  $(X_\alpha)$, $(u_\alpha)$, $(\theta_\alpha)$ 
 such that  $P$  forces  ``$\dot u$ is a $\lambda$-sequence 
of reals, $\dot \theta$ is a $\lambda$-sequence of 
 ordinals, and, for each  $\alpha<\lambda$, $\dot X_\alpha$   
is the  $\dot \theta_\alpha$th element (in the standard 
 order) of  OD($\dot u_\alpha$).''

      Let  $\mu$  be the largest element of  $C$  below  $\lambda$.  
The hypotheses of 
 the theorem imply that  $\mu$  exists and has uncountable cofinality.  (Here 
 and below, we tacitly use Lemma 11 to avoid 
having to say whether cardinals 
 and cofinalities refer to  $V$  or to  $V[G]$; the countable chain condition 
 makes these concepts absolute.)  It follows that, in $V$ where  GCH  holds, 
 $\mu^{\aleph_0}= \mu$.

      We intend to find a set  $M\subseteq\lambda$
 of cardinality  $\mu$,  such that the sets 
 $X_\alpha$   for  $\alpha\in M$  cover ${}^\omega\omega$.  
$M$  will be obtained in the ground model  $V$  as 
 the union of an increasing  $\aleph_1$-sequence of 
approximations  $M_\sigma$   of 
 cardinality  $\leq\mu$  for  $\sigma < \aleph_1$.  
Recall that  $\mu\geq\aleph_1$, so the union  $M$ 
 also has cardinality $\leq\mu$. For limit ordinals  
$\sigma$,  $M_\sigma$   will be the union 
 of the  $M_\tau$'s  for $\tau<\sigma$.  
As  $M_0$   we take the empty set.    The 
 non-trivial part of the construction is the successor step, and for this we 
 need some preliminary work.

      Until further notice, we work in  $V$.  We let $\dot x$ 
range over names for 
 reals, and we identify two names if  $P$  forces that they are equal.  With 
 this convention, the countable chain condition (Lemma 11) allows us to 
 assume that each name $\dot x$ involves only countably many conditions.  
So there is a set  $J\subseteq I$  such that
\par\noindent (a)\quad
all the conditions involved in $\dot x$  have domains  
$F\times n$  with $F\subseteq J$, and
\par\noindent (b)\quad
 $|J| = \aleph_0$.
\par\noindent
 In particular, we can choose such a $J_\alpha$   
for each of the names  $\dot u_\alpha\, (\alpha<\lambda)$ 
 that we chose earlier for the real parameters in ordinal 
definitions of  $\dot X_\alpha$.
 Enlarging  $J_\alpha$, but keeping it countable, 
we can similarly arrange that all 
 conditions involved in  $\dot \theta_\alpha$  
have domains  $F\times n$  with  $F\subseteq J_\alpha$.  
Let  $S$ be the union of these  $\lambda$ countable sets  
$J_\alpha$   and the sets $I_\kappa$   for  
$\kappa\leq \mu$ in  $C$.   So  $|S|=\lambda$.

      Until further notice, consider a fixed but arbitrary set  
$K\subseteq S$  of cardinality  $\mu$  such that  
$I_\kappa\subseteq K$  for all  $\kappa\leq\mu$  in $C$.  
Notice that, 
 for  $\kappa>\mu$  in  $C$, we have  $\kappa > \lambda$ 
(by choice of  $\mu$  and because  $\lambda\notin C$), 
 so  $I_\kappa-S$  and $I_\kappa-K$  have cardinality $\kappa$.
      We shall call $J\subseteq I$  a {\it support} for a name $\dot x$
(of a real) if it 
 satisfies (a) and (b) above and also
\par\noindent (c)\quad
for each  $\kappa\in C$, if $J\cap I_\kappa-K$  is nonempty, 
then it is infinite.
\par\noindent
 Notice that the new clause (c) is easy to satisfy by enlarging $J$.

      Let  $\Cal G$  be the group of those permutations of $I$  that map each $I_\kappa$ into itself and that fix all members of $K$.  
Clearly,  $\Cal G$  acts as a group 
 of automorphisms of the notion of forcing  $P$, by
$$
                         g(p)(g(\kappa,\xi),k) = p(\kappa,\xi,k), 
$$
 and it is well known that such automorphisms also act 
on the class of $P$-forcing names 
(i.e., on the associated Boolean-valued model) and preserve 
 the forcing relation.  It is easy to check that, 
if $J$  supports $\dot x$, then $g(J)$  supports  $g(\dot x)$; 
if, in addition,  $g$  fixes all members of  $J$, 
then it also fixes  $\dot x$.

      If  $J$  is a support then, thanks to clause (c), its 
$\Cal G$-orbit (i.e., 
 its equivalence class under the action of  $\Cal G$) 
is determined by  $J\cap K$  and
$$
\bar J=\{\kappa\in C\mid J\cap I_\kappa-K\neq\emptyset\}.
$$ 
 That is, if  $J'$  is another support with  $J'\cap K = J\cap K $
 and  $\overline{J'}=\bar J$, then 
 there is  $g\in \Cal G$  with  $g(J) = J'$.  
Since  $J\cap K$  is countable and 
 $|K|= \mu = \mu^{\aleph_0}$  there are, 
as  $J$  varies over all supports, only  $\mu$ 
 possibilities for  $J\cap K$.  Also, since  $\bar J$  is 
a countable subset of  $C$ 
 and  $|C|\leq\mu$  (because all uncountable cardinals  
$\leq|C|$  are in  $C$), the 
 number of possibilities for $\bar J$  is  
$\leq\mu^{\aleph _0}=\mu$.  
Therefore, there are only 
 $\mu$  $\Cal G$-orbits of supports.

      For each $\Cal G$-orbit of supports, choose one member 
$J$  such that  $J\cap S =  J\cap K$, 
i.e., such that  $J$  is disjoint from  $S - K$.  Such a  $J$  is easy 
 to find, starting with an arbitrary  $J'$  in the orbit.  
For each  $\kappa\in C$ 
 such that  $J'$  meets  $I_\kappa\cap S - K$, we have  
$\kappa > \lambda$ (for otherwise  $I_\kappa \subseteq K$) 
 and then  $|I_\kappa -S| > \lambda$ so there are 
permutations of  $I_\kappa$ fixing  $I_\kappa\cap K$ 
 pointwise and mapping the (countable) rest of  $J'$ out of  $S$.  Combine such permutations for all relevant  $\kappa$  to get  
$g\in \Cal G$  for which  $J = g(J')$ is as desired.  
Call the  $\mu$  orbit-representatives just chosen 
the {\it standard} supports.

      For any fixed support  $J$, any name $\dot x$  
supported by  $J$  can be specified by giving, for each  $n\in\omega$, a maximal antichain of conditions 
 that are supported by  $J$  and that decide $\dot x(n)$  
and giving those decisions.  
It follows, by  CH, that there are only $\aleph_1$
such names for  each $J$.

      Thus, there are only  $\mu$  names  $\dot x$ 
that have standard supports.  For 
 each of these, fix a countable set  
$A = A(\dot x)\subseteq \lambda$ such that  $P$  forces 
 ``$(\exists\alpha\in\check A)\,\dot x\in \dot X_\alpha$.''  
The existence of such an  $A$  follows from the countable 
 chain condition and the fact that  
``$(\exists\alpha <\lambda)\,\dot x\in \dot X_\alpha$''  is forced.  
Let  $B$  be the union of these sets  $A(\dot x)$  
for all  $\dot x$  with standard support.  As the 
 union of  $\mu$  countable sets,  $B$  has cardinality $\leq\mu$.

      Now un-fix  $K$.  The preceding discussion produces, 
for each  $K\subseteq S $
 of size  $\mu$  with  $I_\kappa\subseteq K$  for  $\kappa\leq\mu$,
a subset  $B$  of $ \lambda$ of cardinality $\leq\mu$.

      At last, we are in a position to complete the definition of the sequence 
 $(M_\sigma)_{\sigma<\aleph_1}$ by carrying out the successor 
step.  Recall that  $\dot u_\alpha$  and  $\dot\theta_\alpha$  are 
 such that  $P$  forces ``$\dot X_\alpha$   is the 
$\dot\theta_\alpha$th set ordinal-definable from  $\dot u_\alpha$''  and 
 that  $J_\alpha$   was chosen so that all conditions 
involved in  $\dot u_\alpha$  and  $\dot\theta_\alpha$  have 
 domains  $F\times n$  with  $F\subseteq J$.  
Now, given  $M_\sigma$, to define  $M_{\sigma+1}$, apply 
 the preceding construction of $B$ from $K$  with
$$
K = K_\sigma=\bigcup_{\alpha\in M_\sigma}J_\alpha
\cup\bigcup_{\kappa\leq\mu\text{ in }C}I_\kappa. 
$$
 As each $J_\alpha$ is countable and $|M_\sigma| \leq\mu$  
we have  $|K| = \mu$  so the 
 construction of $B$  makes sense.  
We set $M_{\sigma+1}=B$  and note that  $|M_{\sigma+1}|\leq\mu$  
as desired.  We also note that the $K_\sigma$ form a 
continuous monotone sequence because the $M_\sigma$ do.

      Having defined $M_\sigma$ for all $\sigma<\aleph_1$
 and thus also their union $M$, we 
 complete the proof of the theorem by showing that, 
for every name  $\dot x$  of a 
 real, $P$ forces ``$(\exists\alpha\in M)\,\dot x\in\dot X_\alpha$''.  
Let  $\dot x$  be given and let $J$  satisfy (a) 
 and (b) in the definition of support for  $\dot x$.  
Let $K_\infty=\bigcup_{\tau<\aleph_1}K_\tau$, so 
 $|K_\infty| = \mu$.  As $J$ is countable, 
we can fix $\sigma<\aleph_1$ such that
$$
J\cap K_\infty \subseteq K_\sigma.\tag1 
$$
 Henceforth, we use notation as in the construction of 
$M_{\sigma+1}$ from $M_\sigma$.  In 
 particular,  $K$  means $K_\sigma$ 
and the notion of support uses this $K$  in clause (c).  
$J$  need not be a support, i.e., clause (c) may fail, but we 
 can enlarge  $J$  to a support by adding elements of  
$I_\kappa-K_\infty$   for all necessary  $\kappa>\mu$.  
(This is possible as, for such  $\kappa$, 
$|I_\kappa|>\lambda>\mu=|K_\infty|$. 
 We needn't worry about  $\kappa\leq\mu$  as 
$I_\kappa\subseteq K_\sigma$  for such  $\kappa$).  
This enlargement preserves \thetag1, so we assume from now on 
that  $J$  is a support.

      Recall that we chose a standard support in every $\Cal G$-orbit of supports. 
 So fix  $g\in\Cal G$  such that  $g(J)$  is standard.  
Neither  $J$  nor  $g(J)$  meets $K_{\sigma+1}-K$.  
In the case of $J$ this follows from \thetag1, while in the case of 
 $g(J)$ it follows from the fact that standard supports don't meet  
$S - K$ (and clearly  $K_\tau\subseteq S$  for all $\tau$).  
Thus, there is  $h\in \Cal G$  such that $h$  agrees 
 with  $g$  on  $J$  and with the identity map on  
$K_{\sigma+1}-K$.  In particular, 
 $h(J) = g(J)$ is standard, and $h$  leaves 
$K_{\sigma+1}$ pointwise fixed (because  all elements of  
$\Cal G$  fix $K$  and $h$  fixes $K_{\sigma+1}- K$.)

      Since  $h(\dot x)$  has standard support  $h(J)$, 
it is one of the  $\mu$  names 
 for which we chose a set  $A = A(h(\dot x))$  
to include in  $B$.  By the defining property of  $A$, $P$ forces 
``$(\exists\alpha\in\check A)\,h(\dot x)\in\dot X_\alpha$ ''  
and thus also
$$
\text{``}(\exists\alpha\in\check A)\,h(\dot x)\text{ is in the } \dot\theta_\alpha\text{th set ordinal-definable from }
\dot u_\alpha.\text{''} 
$$

      For any  $\alpha\in A$, we have  
$\alpha\in B\subseteq M_{\sigma+1}$ 
by definition of  $B$  and  $M_{\sigma+1}$. 
 We also have  $J_\alpha\subseteq K_{\sigma+1}$, 
and so  $h$  fixes $J_\alpha$ pointwise.  It follows, by 
 definition of  $J_\alpha$,  that $h$ fixes the names 
$\dot u_\alpha$  and $\dot\theta_\alpha$.  
So, by \thetag2,  $P$  forces
$$
\text{``}(\exists\alpha\in\check A)\,h(\dot x)\text{ is in the } h(\dot\theta_\alpha)\text{th set ordinal-definable from }
h(\dot u_\alpha),\text{''} 
$$
 and, since  $h$  is an automorphism,
$$
\text{``}(\exists\alpha\in\check A)\,\dot x\text{ is in the } \dot\theta_\alpha\text{th set ordinal-definable from }
\dot u_\alpha.\text{''} 
$$
 Since  $A\subseteq B\subseteq M_{\sigma+1}\subseteq M$, 
we have that  $P$  forces  
``$(\exists\alpha\in\check M)\,\dot x\in \dot X_\alpha$,''  
as required to complete the proof.  
\qed\enddemo

\remark{Remark}
 At the January, 1991, Bar-Ilan conference 
on the set theory  of the reals, 
I described many of the results in this paper and made some 
 conjectures, one of which was that there might be very few uniform 
$\lp1$-characteristics and that 
one might be able to classify them all.  
Shelah promptly informed me that, 
by a countable-support product of forcing 
 notions from [21], he can produce models 
with infinitely many uniform 
$\lp1$-characteristics, 
all of the form  ``the smallest number of $g$-branching 
 subtrees of  ${}^{<\omega}\omega$  needed 
to cover all the paths through an $f$  branching 
 subtree of  ${}^{<\omega}\omega$.'' 
(Here $f$  and $g$  are suitable recursive functions on $\omega$, 
and an $f$-branching tree is one in which each node of level $n$ has 
 exactly $f(n)$ immediate successors.)  These characteristics can be 
 prescribed rather freely, 
and one can get uncountably many of them if one 
 allows boldface $\bp1$, i.e., non-recursive $f$ and $g$.  
I do not know to 
 what extent Shelah's models satisfy 
the additional property, enjoyed by the models in Theorem 9, 
that cardinals not explicitly made to be characteristics
are not even OD$\Bbb R$-characteristics.  
This work of Shelah will appear 
(with some modifications) in [7].
\endremark

      It is clear from Theorem 9 (or even from Theorem 8) 
that we cannot expect restrictive results about  
$\bp2$-characteristics in  ZFC.  But these 
 theorems leave open the possibility of 
restrictive results in stronger theories, perhaps ZFC + MA.  
The following theorem, a corollary of a 
 result of Harrington [8], shows that even Martin's axiom gives no 
 restrictions on the $\boldsymbol\Delta^1_1$-characteristics. 

\proclaim{Theorem 13}
It is consistent with Martin's axiom that $\cc$  be arbitrarily 
 large and that every uncountable cardinal  
$\leq\cc$  be a $\boldsymbol\Delta^1_1$-characteristic and a 
 uniform  $\bS2$-characteristic.
\endproclaim

\demo{Proof}
Harrington [8, Theorem B] obtained models of Martin's axiom in 
 which  $\cc$  is arbitrarily large and there is a  
$\bP2$ well-ordering $\leq$ of a 
 set of reals, having length  $\cc$.  
For any uncountable cardinal  $\kappa<\cc$, let 
 $a$  be the $\kappa$th element in this well-ordering.  
Then the binary relation
$$
R(x,y) \iff y\not\leq a \text{ or } y = x 
$$
is $\bS2$, and a set  $\Cal X$  of reals satisfies  
$\forall y\,(\exists x\in\Cal X)\, R(x,y)$  
if and only if  $\Cal X$  contains all the predecessors of  $a$, 
of which there are  $\kappa$.  This proves 
 the $\bS2$ part of the theorem for  $\kappa < \cc$.
  The case of  $\kappa = \cc$  is trivial 
 as  $\cc$  is a uniform  $\bp1$-characteristic.  
The $\boldsymbol\Delta^1_1$ part follows by 
 Proposition 3(c).  
\qed\enddemo

      Notice that in this proof the well-ordering was used only to 
 produce $\bP2$  sets of arbitrary uncountable cardinality  
$\kappa < \cc$.

\head
5. Duality
\endhead

 There is an intuition that some of the familiar cardinal 
 characteristics of the continuum occur in dual pairs.  
For example, in the abstract of [16], 
Miller refers to dualizing the proof of $\bold{add}(L)\leq\bb$  to 
 obtain  $\dd\leq\bold{cof}(L)$.  
In this section, we make some remarks about this sort 
 of duality, and we attempt to relate it to our theory of 
$\Gamma$-characteristics.

      At first sight, duality seems quite easy to describe.  
Indeed, the  $\ls2$
 relations used in the proof of Proposition 2(b) 
for  $\dd$, $\rr$, ${\bold {cov}}(B)$, and 
 ${\bold {unif}}(L)$  
are precisely the negations of the converses of the $\lp2$ relations 
 used in the proof of Proposition 2(c)  for the dual characteristics  
$\bb$,  $\ss$, ${\bold {unif}}(B)$, and  ${\bold {cov}}(L)$,
 respectively.  Thus, if we define the dual of a binary relation by
$$
\tilde R(x,y) \iff  \neg  R(y,x), 
$$
 then the uniform characteristic determined by $\tilde R$ 
seems to be, in the intuitive sense, 
dual to the uniform characteristic determined by  $R$.

      Some caution is needed, however, with this notion of duality. 
 For example, the $\ls2$-characteristics in Proposition 2(b)  
are in fact uniform $\lp1$-characteristics by Proposition 3(a), 
and in fact the first three of them 
 have quite natural $\lp1$ descriptions.  
(In the definition of  $\dd$, one can 
 replace $\leq^*$ with ``everywhere $\leq$'', 
a similar deletion of ``mod finite'' in 
 the notion of splitting works for $\rr$, 
and for  ${\bold {cov}}(B)$ one can work with 
 nowhere dense closed sets instead of meager sets.)  
Dualizing those definitions leads to a  $\ls1$ form, 
for which the cardinal is $\aleph_0$ 
(and thus not a characteristic, by our definition).  
It would appear that, before dualizing, 
one must be careful to put $R$ into the proper form.  
But what  is the proper form?

      A hint can be obtained by comparing the ``good'', i.e., nicely 
 dualizable definitions of  $\dd$, $\rr$, $\bb$, and $\ss$  
in the proof of Proposition 2 with 
 the ``bad'' versions obtained by deleting ``mod finite''.  
The most evident difference, apart from the 
complexity difference which makes the bad $\lp1$ 
 definitions look better than the good  $\lp2$ ones, 
is that the good  $R$'s  are 
 unchanged by finite modifications of their arguments.  
That is, if  $x$  and $x'$  differ only finitely and if  
$y$  and  $y'$  differ only finitely, then for 
 the good  $R$'s, but not for the bad ones,  
$R(x,y)$  implies  $R(x',y')$.  
We call such  $R$'s invariant.

      What about ${\bold {cov}}(B)$, ${\bold {unif}}(B)$, 
and their measure analogs?  The relations 
 used in the proof of Proposition 2 are not invariant, but they can be 
 replaced with ones that are invariant, have the same complexity, 
and lead to the same duals.

      For the category situation, 
we define a new coding of meager sets by 
$$
 x\in \Meager'(y) \iff  
(\exists s,t\in {}^{<\omega}\omega) \,x * s\in \Meager(y*t), 
$$
where \Meager\ is as in the proof of Proposition 2(b) 
and where  $x * s$  is  $x$  but with $ s$  in 
 place of the initial segment of the same length in  $x$, i.e.,
$$
(x*s)(n) =\cases
s(n)&\text{if }n\in\text{domain}(s)\\
x(n)&\text{otherwise.}
\endcases
$$
 Thus, the relation  ``$x\in\Meager'(y)$''  is  
``$x\in \Meager (y)$'' enlarged just enough to be invariant.  
Notice that $\Meager'(y)$  is meager 
(since there are only countably many possibilities for $s$ and $t$) 
and every meager set has a superset of the form  $\Meager'(y)$.  
Furthermore,  ``$x\in \Meager'(y)$''  is a
$\ls2$  relation, so $\Meager'$ 
can replace $\Meager$ in the proof of Proposition 2.

      For the measure situation, things are a bit more complicated. 
 Defining $\Null'$ in exact analogy with $\Meager'$ 
would make ``$x\in \Null'(y)$'' 
 $\ls3$  because  ``$x\in \Null(y)$''  is $\lp2$, 
so the complexity would increase. 
 Changing  $(\exists s,t)$ to $(\forall s,t)$  
would make $\Null'(y)$  empty, since there is 
 always a  $t$  for which  $\Null(y*t)$  is empty.  
We observe however, that we can safely use  $\forall s$; 
that is, if we define
$$
x\in \Null^* (y) \iff  
(\forall s\in{}^{<\omega}\omega)\,x*s\in \Null(y), 
$$
 then this relation is $\lp2$, invariant with respect to  $x$, 
and still enjoys the crucial property that every set  $A$  
of measure zero is included in one 
 of the form  $\Null^* (y)$.  (For the proof, simply observe that 
 $\{x*s\mid x\in A,s\in{}^{<\omega}\omega\}$  
has measure zero and is therefore included in some  $\Null(y)$.)
  To obtain invariance with respect to  $y$, we use a different 
 coding.  Every set  $A$  of measure zero can be covered 
by a sequence of sets $A_n$ each of which is a 
union of finitely many basic open sets  $B_s$  (as in 
 the proof of Proposition 2) and has measure below 
some prescribed positive  bound $\varepsilon_n$.  
Dovetailing infinitely many such constructions with suitable 
$\varepsilon$'s, we can find a sequence of sets  
$A_n$   such that  $A_n$   is a finite union 
 of basic open sets,  $A_n$   has measure $\leq2^{-n}$, 
and every element of  $A$  is in 
 infinitely many  $A_n$'s.  (This construction is in [2, Lemma 1.1].) 
 Conversely, if each $A_n$   has measure  $\leq2^{-n}$, 
then $\{x\mid x\in A_n\text{ for infinitely many }n\}$ 
has measure zero.  
Let  $f:\omega\times\omega\to[{}^{<\omega}\omega]^{<\omega}$
 be a recursive function such that, 
if $n$ is fixed and $k$ varies over $\omega$, $f(n,k)$ 
 enumerates all finite sets  $F\subseteq{}^{<\omega}\omega$  
such that $\bigcup_{s\in F}B_s$ has measure $\leq2^{-n}$. 
 Then 
$$
x\in \Null^\dag (y) \iff  
(\forall m)\,(\exists n\geq m)\,
x\text{ has a initial segment in }f(n,y(n)) 
$$
 defines a $\lp2$ relation, invariant with respect to  $y$, 
such that all sets of the form  $\Null^\dag (y)$  
have measure zero and all sets of measure zero have 
 supersets of this form.  Finally, we obtain the desired $\Null'$, invariant in both variables, by putting  $\Null^\dag$  
 in place of  $\Null$  in the defintion 
 of  $\Null^*$:
$$
x\in \Null'(y) \iff  
(\forall s\in{}^{<\omega}\omega)\,x*s\in \Null^\dag(y), 
$$

      The preceding discussion suggests that a suitable context 
for duality  is uniform characteristics given by invariant relations.  
We leave it to  the reader to verify that the relation 
$\Meager'(y)\subseteq\Meager'(x)$  and its 
 dual  $\Meager'(x)\not\subseteq\Meager'(y)$  determine  
$\bold{cof}(B)$  and $\bold{add}(B)$,  respectively, 
so these are dual characteristics determined by invariant 
 relations.  
The analog for measure also holds, as does the improvement, 
 described in Section 2, from $\lS1$ to $\lp2$ for $\bold{add}(L)$
  and from $\lP1$  to $\ls2$  for $\bold{cof}(L)$, 
because the relations obtained from Theorems 0.9 and 0.10 of 
 [1] are invariant.

      We close this section with a very easy result, dual to Theorem 5, 
as  propaganda for this notion of duality.
\proclaim{Proposition 14}
Let $\kappa$ be the uniform $\bp2$-characteristic 
determined by an  invariant $\bp2$ relation  $R$.  
Then $\kappa\leq{\bold {unif}}(B)$.
\endproclaim
\demo{Proof}
Since  $R$  determines a characteristic, there must be, for each 
 $y\in{}^\omega\omega$, at least one  $x\in  {}^\omega\omega $
 such that  $R(x,y)$.  But then, being 
 invariant under finite modifications,
$$
R_y  = \{x\in {}^\omega\omega\mid R(x,y)\}
$$ 
 is dense.  It is also $\bp2$, i.e., a  $G_\delta$   set, 
so it is comeager.  Thus, letting $\Cal X$ be a non-meager set
of the smallest possible size ${\bold {unif}}(B)$, we have that 
$\Cal X$ meets $R_y$ for each $y$.  
Thus, by definition of  $\kappa$, we have 
 $\kappa\leq|\Cal X| = {\bold {unif}}(B)$.
\qed\enddemo

      In particular, we have the following analog of Corollary 6,
 including the other half of Rothberger's theorem [20] 
along with some easier known results.

\proclaim{Corollary 15}
All of  $\bb$, $\ss$, and  ${\bold {cov}}(L)$  are 
$\leq{\bold {unif}}(B)$. \qed
\endproclaim

      We remark that, in the proof of Proposition 14, 
invariance of  $R$  was 
 used only with respect to the variable  $x$.

\head
6. Partition Characteristics
\endhead

      This section is devoted to some characteristics 
connected with partition  theorems.  
Some of these characteristics first attracted my attention as 
 possible new examples of uniform $\lp1$-characteristics.  
(This was before  Shelah showed how to get a plentiful supply 
of uniform $\lp1$-characteristics [7].)  
Others arose as duals.  
They seem to have some intrinsic interest, 
so we present here what is known about them.
      
We begin with Ramsey's theorem in the simple form:  
If  $[\omega]^2$   is  partitioned into two pieces, then 
there is an infinite  $H\subseteq\omega$  such that 
$[H]^2$  is included in one piece.  
As usual,  $[X]^n$   is the set of $n$-element 
 subsets of $X$, and an $H$  as in the theorem is said to be 
{\it homogeneous} for the partition.  
We call $H$  {\it almost homogeneous} if there is a 
finite  $F\subseteq H$ such that  $H - F$  is homogeneous.
      
Let
$$\align
 {\bold {par}} =&\text{minimum number of partitions }
 \Pi: [\omega]^2  \to 2\text{ such that}\\
&\text{ no single } H\in[\omega]^\omega\text{ is almost 
homogeneous for them all}\\
 \intertext{and}
 {\bold {hom}} =&\text{minimum number of infinite subsets }
H \text{ of } \omega\\
&\text{ such that every 
           partition } \Pi: [\omega]^2  \to 2\\
&\text{ has an almost homogeneous set among 
           these } H\text{'s}.
\endalign$$
 Both of these are easily seen to be uncountable.  
Since ``$H$  is almost 
 homogeneous for $\Pi$'' is an invariant $\ls2$  relation, 
${\bold {par}}$ is a uniform  $\lp2$-characteristic,
 ${\bold {hom}}$ is a uniform $\ls2$-\ 
(hence $\lp1$-)characteristic, and their 
 definitions are dual to each other. (That ${\bold {hom}}$ is 
a uniform  $\lp1$-characteristic can also be seen by observing 
that it is unchanged if ``almost'' is removed from the definition; 
the same change would turn ${\bold {par}}$ into $\aleph_0$.)

      Our first result is that ${\bold {par}}$ is nothing new.

\proclaim{Theorem 16}
${\bold {par}}$ is the smaller of  $\bb$  and  $\ss$.
\endproclaim
 
\demo{Proof}
 First, we consider any $\kappa < {\bold {par}}$, 
i.e., any  $\kappa$  such that every  $\kappa$ partitions  
$[\omega]^2\to2$  
have a common, infinite, almost homogeneous set.

We claim that  $\kappa<\bb$.  To prove this, 
let a family  $\Cal F$ of  $\kappa$ non-decreasing functions 
$f:\omega\to\omega$  be given; we seek a single $g$ eventually 
 majorizing them.  Each  $f\in \Cal F$ induces a partition
 $\Pi_f: [\omega]^2\to 2$, namely
$$
\Pi_f\{a < b\} = 0 \iff  f(a) < b. 
$$
 A homogeneous set of color 1 for $\Pi_f$ must be finite, 
being bounded by  $f$  of its first element.  
So the common, infinite, almost homogeneous set $H$ for 
 all the $\Pi_f$, $f\in\Cal F$,  must be 
almost homogeneous of color  0.  That is, for 
 each  $f\in\Cal F$,  we have  $f(a) < b$  
for all sufficiently large  $a < b$  in  $H$. 
 It follows that the function sending each  $n\in \omega$  
to the second element of  $A$ after  $n$  
eventually majorizes each  $f\in\Cal F$.

      We claim further that  $\kappa < \ss$.  
Let a family $\Cal S$  of  $\kappa$ infinite subsets 
$S$ of $\omega$  be given; we seek an infinite set 
not split by any of them.  Each $S\in\Cal S$  
induces a partition $\Pi_S: [\omega]^2  \to 2$, namely
$$
\Pi_S\{a < b\} = 0  \iff   a\in S. 
$$
 Clearly, a set almost homogeneous for $\Pi_S$ is 
not split by  $S$, so the 
 hypothesis on  $\kappa$  provides the desired unsplit set.

The preceding two claims establish that  
${\bold {par}} \leq \min\{\bb,\ss\}$.

      To prove the converse, consider any  $\kappa < \min\{\bb,\ss\}$,
 and let a family of  $\kappa$  partitions 
$\Pi_\alpha: [\omega]^2\to2$, for $\alpha<\kappa$, be given.  
We seek a set almost homogeneous for all the $\Pi_\alpha$'s.  
For each $\alpha< \kappa$  and each  $a\in\omega$, let
$$
S_{\alpha,a}   = \{b\in\omega-\{a\}\mid\Pi_\alpha\{a,b\} = 0\}. 
$$
 Since there are only  $\kappa\cdot\aleph_0<\ss$  sets  
$S_{\alpha,a}$, they do not form a splitting  family.  
So let  $A$  be an infinite set not split by any of them.  This means 
 that, for each $\alpha$ and $a$, the value of 
$\Pi_\alpha\{a,b\}$  is the same, say  $v_\alpha(a)$, 
 for all sufficiently large  $b\in A$,  
say all such that  $b > f_\alpha (a)$.

      The same argument, applied within $A$ to the  
$\kappa$ sets  $\{a\mid v_\alpha (a) = 0\}$, 
 provides an infinite  $B\subseteq A$  such that  
$v_\alpha(a)$  has the same value, say  $i_\alpha$, 
 for all sufficiently large $a\in B$, say all such  $a\geq u_\alpha$.
Since  $\kappa < \bb$, the  $\kappa$  functions $f_\alpha$   
are all eventually majorized by a single function  $g$. 
 Increasing  $u_\alpha$ if necessary, we can arrange that  
$g(a)\geq f_\alpha(a)$  for all  $a\geq u_\alpha$.  
Finally, we construct the desired almost homogeneous infinite set 
 $H\subseteq B$  by choosing its members inductively from  
$B$  so that, if  $a < b$  are 
 in  $H$, then  $g(a) < b$.  To see that  $H$  is 
almost homogeneous for each $\Pi_\alpha$, 
 suppose  $a$  and  $b$  are in  $H$  and $u_\alpha\leq a < b$.  
Then  $f_\alpha (a)\leq g(a) < b$  and 
 so  $\Pi_\alpha\{a,b\} = v_\alpha (a) = i_\alpha$.   
\qed\enddemo

     Attempting to dualize the preceding argument, 
to obtain  ${\bold {hom}} = \max\{\dd,\rr\}$, 
 we succeed only partially.  To state the result that we obtain, 
we need the following variant of  $\rr$  
introduced and studied by Vojt\'a\v s [24,25].
$$\align
\rr_\sigma  =&\text{ smallest cardinality of any\ }
  \Cal X\subseteq[\omega]^\omega\\
&\text{ such that, for any countably many sets\ }
 Y_n \in [\omega]^\omega,\ n\in \omega,\\
&\text{ there is  }X\in\Cal X\text{  not split by any }Y_n.
\endalign
$$
 It is clear that $\rr_\sigma$ is a uniform $\lp3$-characteristic 
and that  $\rr_\sigma\geq \rr$.  It 
 is an open problem whether $\rr_\sigma=\rr$ (provably in ZFC).  
$\rr_\sigma$ arises naturally  in analysis as the characteristic 
associated to the Bolzano-Weierstrass theorem; 
it is the smallest cardinality of any  
$\Cal X\subseteq [\omega]^\omega$   such that every  bounded 
sequence $(x_n)_{n\in\omega}$ of real numbers has a convergent 
subsequence of  the form $(x_n)_{n\in X}$ with $X\in\Cal X$ [25].

      Observe that the $\lp3$ relation determining  $\rr_\sigma$, 
namely ``no term of the  sequence coded by  $y$ splits  $x$,''  
is invariant (for reasonable coding), so 
 it makes sense to consider the dual relation, 
defining the uniform $\ls3$-\  (hence  $\lp2$-)characteristic 
``the minimum number of $\omega$-sequences of sets such 
 that every infinite set is split by some term of 
one of these $\omega$-sequences.'' 
 But this is simply $\ss$.  Thus,  $\ss$  is dual to both  
$\rr$  and $\rr_\sigma$, a  circumstance that helps 
to explain why the following dual of Theorem 16 
 looks weaker than one might expect.

\proclaim{Theorem 17}
$\max\{\rr,\dd\}\leq{\bold {hom}}\leq\max\{\rr_\sigma,\dd\}$.
\endproclaim

\demo{Proof}
For the first inequality, fix a family $\Cal X$  of 
${\bold {hom}}$ infinite subsets  $H$  of  $\omega$, 
containing an almost homogeneous set for every 
$\Pi:[\omega]^\omega\to 2$. 
 For any $S\subseteq\omega$,  define $\Pi_S$   
as in the proof of Theorem 16, and observe 
 that an  $H$  almost homogeneous for $\Pi_S$ is not split by  $S$. 
 Thus,  $\Cal X$  is unsplittable and so $\rr\leq{\bold {hom}}$.  
Similarly, given  a non-decreasing  $f\in{}^\omega\omega$, 
 define $\Pi_f$  as in the proof of Theorem 16, 
and observe that, if $H$ is  homogeneous for $\Pi_f$   
then, as in that proof,
$$
g_H(n) = \text{the second element of } H \text{ after } n 
$$
 defines a  $g_H$   eventually majorizing  $f$.  
Thus  $\{g_H\mid H\in\Cal X\}$ is a dominating 
 family and so  $\dd\leq{\bold {hom}}$.

      For the second inequality, let  
$\kappa = \max\{\rr_\sigma,\dd\}$.  
Let  $\Cal X$  be a family of $\rr_\sigma$   sets as in the definition 
of  $\rr_\sigma$.
  Inside each  $X\in\Cal X$, let $\Cal Y(X)$  be an 
 unsplittable family of $\rr$ sets.  
Let $\Cal Y = \bigcup_{X\in\Cal X}\Cal Y(X)$,  
and let $\Cal D\subseteq{}^\omega\omega$  be a 
 dominating family of cardinality $\dd$.  
For each  $Y\in\Cal Y$ and $f\in \Cal D$,  let $Z = Z(Y,f)$  be an infinite subset of  $Y$  such that, if  $a < b$  are in  $Z$ 
 then  $f(a) < b$.  Since all of  $\rr$, $\rr_\sigma$, and $\dd$  
are  $\leq\kappa$, there are at most  $\kappa$ 
 sets  $Z(Y,f)$.  
We shall complete the proof by showing that every partition
$\Pi:[\omega]^\omega\to 2$  has an almost homogeneous set 
among the  $Z(Y,f)$'s.

      So let $\Pi$ be given, and consider the countably many sets
$$
S_a  = \{b\in\omega-\{a\}\mid\Pi\{a,b\} = 0\} 
$$
 for  $a\in\omega$.  By choice of $\Cal X$,  find $X\in\Cal X$  
not split by any  $S_a$.  Thus, for each  $a$,  
there are $v(a)\in \{0,1\}$  and  $g(a)\in\omega$ 
 such that, whenever  $b\in X$  and  $b\geq g(a)$,  
then $\Pi\{a,b\}=v(a)$.  By choice of $\Cal Y(X)$,  find 
 $Y\in\Cal Y(X)\subseteq\Cal Y$  not split by $\{a\mid v(a) = 0\}$.
 Thus, there are  $i\in \{0,1\}$  and 
 $u\in\omega$  such that, if $a\in Y$  and  $a\geq u$,  
then  $v(a) = i$.  By choice of  $\Cal D$, 
 find  $f\in \Cal D$ eventually majorizing  $g$; 
increasing  $u$  if necessary, we may 
 assume that  $f(a) > g(a)$  for all  $a\geq u$.  
Now if  $a$  and  $b$  are in  $Z(Y,f)$ 
 and  $u < a < b$, then, $g(a) < f(a) < b$  (by definition of  $Z(Y,f)$) and 
 therefore $\Pi\{a,b\} = v(a) = i$.  
So  $Z(Y,f)$  is almost homogeneous for $\Pi$, as required. 
\qed\enddemo
 
 We generalize ${\bold {par}}$ and ${\bold {hom}}$  
by considering partitions of  $[\omega]^k$   instead  of  $[\omega]^2$.  
(One could also consider partitions into a larger (finite) number 
 of pieces, but it is easy to check that this would not affect either 
 characteristic.)  Let  ${\bold {par}}_k$   and  ${\bold {hom}}_k$   
be defined exactly like ${\bold {par}}$ and ${\bold {hom}}$
 except that $[\omega]^2$   is replaced by  $[\omega]^k$.  
Notice that  ${\bold {par}}_1=\ss$  and ${\bold {hom}}_1=\rr$.  
Henceforth, we consider only  $k\geq2$.  The proofs of Theorems 16 
 and 17 generalize easily to these higher values of  $k$, 
but in fact one can say slightly more, as was pointed out
 to me by Laflamme who attributed the observation to Shelah.

\proclaim{Proposition 18}
 ${\bold {par}}_k=\min\{\bb,\ss\}$  and 
${\bold {hom}}_k=\max\{\rr_\sigma,\dd\}$ for $k\geq3$.
\endproclaim
 
\demo{Proof}
In view of the preceding remarks and the obvious fact that 
${\bold {hom}}_k\leq{\bold {hom}}_l$   for  $k\leq l$,  
all we need to prove is that $\rr_\sigma\leq{\bold {hom}}_3$.  
Let $\Cal X $
 be a family of ${\bold {hom}}_3$ infinite subsets of $\omega$,
 containing almost homogeneous 
 sets for all partitions $\Pi: [\omega]^3\to2$.  
We claim that $\Cal X$  is as required in 
 the definition of $\rr_\sigma$.  
Let countably many sets  $Y_n$   be given.  Define
 $\Pi : [\omega]^3\to2$  by
$$
\Pi\{a < b < c\} = 0  \iff(\forall n\leq a)\, (b\in Y_n\iff  c\in Y_n),
$$ 
 and let  $H\in\Cal X$  be almost homogeneous for $\Pi$; 
deleting finitely many elements from  $H$  we get a 
homogeneous set  $H'$, and we complete the proof by 
 showing that  $H'$ (and therefore also  $H$) 
is not split by any  $Y_n$.  
Let  $a$ be the smallest element of  $H'$.

      If $\Pi$  maps $[H']^3$   to  1, then infinitely many sets
$$
\{n \leq a\mid b\in Y_n \} 
$$
 are all distinct, as  $b$  varies over  $H' - \{a\}$,  
but they are all subsets of 
 the finite set  $\{0,1,...,a\}$, so this is absurd.  
Therefore  $\Pi$  maps $[H']^3$  to  0.  
This means that, for any  $n$,  $Y_n$   contains all or none of those 
 $b\in H'$ that are greater than the next element of  $H'$ after  $n$.
  So $Y_n$   does not split $H'$.  
\qed\enddemo
 
      Going beyond Ramsey's theorem, we can define 
analogous characteristics  associated with the partition theorems 
of Nash-Williams [19], Galvin and  Prikry [6], and Silver [22].  
Little is known about these characteristics, 
 but we list for reference some elementary facts.  
As we go from weaker to  stronger partition theorems, 
the ${\bold {hom}}$ characteristics weakly increase and in 
 particular are all $\geq\max\{\rr_\sigma,\dd\}$,  and the ${\bold {par}}$ 
characteristics weakly 
 decrease and in particular are all $\leq\min\{\bb,\ss\}$.  
A lower bound for the ${\bold {par}}$
 characteristics is the distributivity number $\bold h$ 
defined as follows.  Call a family  
$\Cal D\subseteq[\omega]^\omega$  {\it dense open} if  
(a) every  $X\in[\omega]^\omega$ has a subset in $\Cal D$ and  
(b) if  $X\subseteq^*Y$  and  $Y\in \Cal D$  then  $X\in \Cal D$.  
Then  $\bold h$  is the smallest cardinal  $\kappa$  such that 
some  $\kappa$  dense open families have empty intersection. 
 That $\bold h\leq$ the analogs of ${\bold {par}}$ 
associated to various partition theorems  follows easily 
from the fact that those theorems ensure that, for any 
 partition of the appropriate sort, the almost homogeneous sets form a dense  open family.  
(Duality provides an upper bound for the characteristics 
 analogous to ${\bold {hom}}$, namely the smallest cardinality 
of a family $\Cal X$  of sets  that meets every dense family.  
Unfortunately, this cardinal equals  $\cc$, a 
 trivial upper bound.)

      It is clear that one can similarly associate characteristics 
with other  partition theorems, for example 
the canonical partition theorem of Erd\H os and Rado [4] 
or the finite sum theorem of Hindman [10].  We shall discuss only 
 two more analogs of ${\bold {par}}$ and two analogs of 
${\bold {hom}}$, associated to very weak partition theorems.

      The first of these theorems can be viewed as the 
``canonical partition  theorem for singletons,'' so we denote the 
analogs of ${\bold {par}}$ and ${\bold {hom}}$ with the 
 subscript  $1c$, but it's really just the (infinitary) pigeonhole principle: 
 If  $f : \omega\to\omega$, then there is an infinite  
$H\subseteq\omega$  on which  $f$  is constant 
 or one-to-one.  We define
$$\align      
{\bold {par}}_{1c}=&\text{smallest cardinality of any  }
\Cal X\subseteq{}^\omega\omega\text{  such that there is no}\\
& H\in[\omega]^\omega\text{ such that each  }f\in \Cal X \\ 
&\text{is constant or one-to-one
                on a cofinite subset of } H,\\
 \intertext{and}
{\bold {hom}}_{1c}=&\text{smallest cardinality of any }
\Cal X\subseteq[\omega]^\omega\text{ such that every }\\
&f\in {}^\omega\omega\text{ is constant or one-to-one on some } 
H\in\Cal X.
\endalign
$$
 As usual, ${\bold {hom}}_{1c}$ would be unaffected if we included 
``mod finite'', and then it is clearly dual to ${\bold {par}}_{1c}$.
 
\proclaim{Proposition 19}
(a)\quad  ${\bold {par}}_{1c}= \min\{\bb,\ss\}$.\newline
(b)\quad  $\max\{\rr,\dd\}\leq{\bold {hom}}_{1c}\leq
\max\{\rr_\sigma,\dd\}$.
\endproclaim

\demo{Proof}
Notice that  $f$  is one-to-one or constant on $H$  if and only if 
 $H$  is homogeneous for the partition of $[\omega]^2$   
that sends $\{a,b\}$ to  0  if and only if  $f(a) = f(b)$.  
This immediately implies  ${\bold {par}}_{1c}\geq{\bold {par}}$  and 
${\bold {hom}}_{1c}\leq{\bold {hom}}$.  In view of Theorems 16 
and 17, we have half of each of (a) and  (b).  
It remains to prove  ${\bold {par}}_{1c}\leq\bb,\ss$  and  
${\bold {hom}}_{1c}\geq\rr,\dd$.  The parts 
 pertaining to $\ss$  and $\rr$  follow from the observation that 
an infinite set $H$  is not split by  $X$  
if and only if the characteristic function of  $X$  is 
 constant or one-to-one on a cofinite subset of  $H$. 
(It can't be one-to-one.)

      To prove  ${\bold {par}}_{1c}\leq\bb$, 
consider an arbitrary  $\kappa < {\bold {par}}_{1c}$, and let a family 
$\Cal F$ of $\kappa$  functions  $f\in{}^\omega\omega$  be given; 
we must find a single  $g\in{}^\omega\omega$
 eventually majorizing them all.  For each  $f\in\Cal F$,  
partition $\omega$ into finite  intervals  $[a_0,a_1)$, $[a_1,a_2)$, 
etc., where  $0=a_0<a_1<a_2<\dots$   and  $a_{n+1}>f(n)$.  
Define  $\hat f\in{}^\omega\omega$ by letting  
$\hat f(k)=n$  for all $k\in [a_n,a_{n+1})$.
 As  $\kappa < {\bold {par}}_{1c}$, find an infinite  
$H\subseteq\omega$  such that each $\hat f$  is constant or 
 one-to-one on a cofinite subset of $H$.  
Define  $g\in{}^\omega\omega$  by letting  $g(n)$  be 
 the $2^n$th element of  $H$; 
we shall show that  $g$  eventually majorizes every 
 $f\in\Cal F$.  Fix any  $f\in\Cal F$.  
As  $f$  is not constant on any infinite set, the 
 defining property of  $H$  ensures that only finitely many 
of the intervals  $[a_n,a_{n+1})$  associated to  
$f$  meet  $H$  more than once.  It follows that, for  sufficiently 
large $n$,  $g(n)$  is in an interval later than the $n$th, so 
 $g(n)\geq a_{n+1}>f(n)$.  
This completes the proof that  $\kappa <\bb$  and therefore 
 ${\bold {par}}_{1c}\leq\bb$.

      The proof that ${\bold {hom}}_{1c}\geq\dd$  is quite similar.  
Let $\Cal X$  be as in the  definition of ${\bold {hom}}_{1c}$, 
and associate to each  $H\in\Cal X$  the function  $g$ 
 defined as above, sending  $n$  to the $2^n$th element of  $H$. 
 To see that these 
${\bold {hom}}_{1c}$ functions $g$ form a dominating family, 
consider any  $f\in{}^\omega\omega$, define $\hat f$ 
 as above, find  $H\in\Cal X$  such that $\hat f$  is 
one-to-one or constant on a cofinite 
 subset of $H$, and argue as above that $g$ eventually majorizes $f$. 
\qed\enddemo

The last pair of partition characteristics that we discuss is defined 
 like the pair  ${\bold {par}}_{1c}$   and  ${\bold {hom}}_{1c}$
  except that ``one-to-one'' is weakened to ``finite-to-one''.  
We call these  ${\bold {par}}_{1cf}$ and  ${\bold {hom}}_{1cf}$,
 where  f  stands for  ``finite''.  
Clearly  ${\bold {par}}_{1cf}\geq{\bold {par}}_{1c}$ and  
${\bold {hom}}_{1cf}\leq{\bold {hom}}_{1c}$.

\proclaim{Proposition 20}
(a)\quad ${\bold {par}}_{1cf}= \ss$.\newline
(b)\quad $\rr\leq{\bold {hom}}_{1cf}\leq\rr_\sigma$.
\endproclaim

\demo{Proof}
The proofs that ${\bold {par}}_{1cf}\leq\ss$  and  
${\bold {hom}}_{1cf}\geq\rr$  are the same as the 
corresponding proofs for ${\bold {par}}_{1c}$ and 
${\bold {hom}}_{1c}$; just replace ``one-to-one'' with 
 ``finite-to-one''.

To see that ${\bold {hom}}_{1cf}\leq\rr_\sigma$, we let 
$\Cal X$  be as in the definition of $\rr_\sigma$ and we show that 
$\Cal X$  also has the property required in the definition of 
${\bold {hom}}_{1cf}$.  Let $f\in{}^\omega\omega$,  and let  
$Y_n=f^{-1}\{n\}$  for each $n\in\omega$.  By hypothesis, 
$\Cal X$  contains an infinite set $H$ not split by any $Y_n$.  
If, for some $n$, $H$ is almost included in $Y_n$, 
then $f$ is constant with value  $n$ on a cofinite part of  $H$.  
Otherwise,  $H$  is almost disjoint from every  $Y_n$, 
 and this means that  $f$  is finite-to-one on  $H$.

      The proof that ${\bold {par}}_{1cf}\geq\ss$  is quite similar.  
We consider any  $\kappa < \ss$ and show that 
$\kappa<{\bold {par}}_{1cf}$.  Let  $\kappa$  functions  
$f\in{}^\omega\omega$ be given.  The 
 $\kappa\cdot\aleph_0<\ss$  sets $f^{-1}\{n\}$, for the given $f$'s 
and all $n\in\omega$, do not form a splitting family, 
so let  $H$  be an infinite set not split by any of them. 
 For each of the given $f$'s, the argument in the preceding paragraph shows  that  $f$  is finite-to-one or constant 
on a cofinite subset of  $H$. 
\qed\enddemo

\head
7. Questions
\endhead
 
 1.\quad   Among the familiar cardinal characteristics of the
continuum [23], the distributivity number $\bold h$ 
and the closely related groupwise density number  $\bold g$  
do not seem to fit our definition of $\Gamma$-characteristics, 
because their definitions involve counting (dense or 
groupwise dense) families of reals rather than counting reals.  
Can one give equivalent  definitions of $\bold h$  and $\bold g$
  showing that they are (at least) OD$\Bbb R$-characteristics?  
(Of course a smaller $\Gamma$ than OD$\Bbb R$ would be 
  preferable.)

 2.\quad  Many more of the familiar characteristics are 
$\Gamma$-characteristics for a reasonable  $\Gamma$  
but are not known to be uniform $\Gamma$-characteristics for any 
$\Gamma$.  Examples include $\bold p$, $\ttt$, $\aa$, $\ii$, 
and $\bold u$.  Are any of these provably uniform 
OD$\Bbb R$-characteristics?

 3.\quad   We saw in Section 2 that $\bold{add}(B)$  and 
$\bold{cof}(B)$ are uniform $\lp2$-characteristics.  
Duality suggests that it should be possible to 
replace $\lp2$ with $\ls2$ and therefore with $\lp1$ 
for one of the two. Theorem 5 requires that one to be 
$\bold{cof}(B)$, since $\bold{add}(B)$ can 
consistently be $<{\bold {cov}}(B)$.  Therefore, we ask:  
Is $\bold{cof}(B)$ a uniform $\lp1$-characteristic?

 4. \quad  To what extent are the hypotheses about  $C$  in 
Theorem 9 needed for the theorem and not just for our proof?  
We remarked before stating the theorem that 
$C$ has to be closed under limits of cofinality $\omega$.  
If, as in our proof, each  $\kappa\in C$  is the cardinality of a 
maximal almost disjoint family, then, by a result of 
Hechler [9, Thm 1], $C$  must be closed under singular limits.  
But there might be proofs that don't rely on 
maximal almost-disjoint families and allow non-closed sets  $C$. 
The requirement that  $C$  contain the immediate successors of all its members of cofinality  $\omega$  
cannot be deleted entirely, as  $\max(C)$, which is to be  
$\cc$  in the extension, had better not have cofinality $\omega$. 
But one might be able to significantly weaken it.  
And the requirement that  $C$  contain all uncountable cardinals  
$\leq|C|$  is purely a technical requirement for our proof.

 5.\quad   Is  $\bold{hom} = \bold{hom}_3$?  One way to settle this (affirmatively) would be to prove $\rr=\rr_\sigma$, 
but  $\bold{hom} = \bold{hom}_3$ might be easier to prove.  
(A meta-question:  Clarify the connection between this question
 and the ``2 versus 3'' problem in the theory of initial segments of 
models of arithmetic [13, p. 226].)  
A related question, bringing the $\rr$  versus $\rr_\sigma$ question 
to the forefront without the extra complication of $\max\{-,\dd\}$, 
is whether either of the inequalities in Proposition 20(b) is 
reversible.

The referee has pointed out a similar problem concerning the
 cardinal ${\bold {cov}}(L)$.  Define ${\bold {cov}}_\sigma(L)$  
to be the minimum cardinality for a family $\Cal X$  of 
 measure zero sets of reals such that every countable set of reals is 
a subset of some $X\in\Cal X$ .  
Is  ${\bold {cov}}_\sigma(L)$ equal to ${\bold {cov}}(L)$?

\head
Appendix. Shelah's Proof of $\dd\leq \ii$
\endhead

An infinite family  $\Cal I\subseteq[\omega]^\omega$ 
is said to be independent if, whenever 
 $\Cal X$  and  $\Cal Y$  are disjoint, finite subfamilies of $\Cal I$, 
then the intersection
$$
 W(\Cal X,\Cal Y) = (\bigcap_{X\in\Cal X}X)
\cap(\bigcap_{Y\in\Cal Y}(\omega-Y)) 
$$
 is infinite.  (Most authors only require that  $W(\Cal X,\Cal Y)$  
be nonempty, but when $\Cal I$ is infinite this definition is 
equivalent to ours, and we don't wish to
 consider finite independent families.)  
The characteristic $\ii$  is defined to 
 be the smallest cardinality of any maximal independent family.

\proclaim{Theorem 21}
{\rm (Shelah, [23, appendix])} $\dd\leq\ii$.
\endproclaim

The following proof is based on the one in [23], but it avoids a few of 
 the complications in that proof.  
Claude Laflamme has informed me that Bill Weiss has simplified
Shelah's proof in a very similar way.
We begin with a lemma that is 
essentially  due to Ketonen [12, Prop. 1.3].

\proclaim{Lemma 22}
Let $C_n$   be a decreasing sequence of infinite subsets of 
$\omega$, and let  $\Cal A$ be a family of fewer than 
$\dd$  subsets of $\omega$ such that each $A\in \Cal A$ 
 has infinite intersection with each $C_n$.  Then there is 
a subset $B$  of $\omega$ such that $B\subseteq^*C_n$
for all $n$  and $A\cap B$  is infinite for all  $A\in \Cal A$.
\endproclaim

\demo{Proof}
For any  $h:\omega \to \omega$, let
$$
B_h=\bigcup_{n\in\omega}(C_n \cap h(n)). 
$$
 As the  $C_n$ form a decreasing sequence, it is clear that 
$B_h \subseteq^*C_n $  for all $n$.  
Our goal is to choose $h$  so that  $A\cap B$   is infinite 
for all  $A\in \Cal A$. 
 Define, for  $A\in\Cal A$ and  $n\in \omega$,
$$
f_A(n) = \text{the }n\text{th element of }A\cap C_n. 
$$
 Notice that, if $h(n) > f_A(n)$  for a particular $A$ and $n$, then 
$|A\cap B_h|\geq n$, as 
$A\cap B_h\supseteq A\cap C_n \cap (f_A(n) + 1)$.  
So it suffices to choose $h$ so that, for each $A\in\Cal A$, 
infinitely many  $n$  satisfy  $h(n)\geq f_A(n)$, 
 i.e., $h\not\leq^*f_A$.  As  $|\Cal A|<\dd$, the functions $f_A$ for
$A\in \Cal A$  cannot constitute a dominating family, 
so such an  $h$  exists.
\qed\enddemo

\demo{Proof of Theorem 21}
Suppose  $\Cal I$ is an independent family of cardinality 
 smaller than $\dd$.  We shall show that  
$\Cal I$  is not a maximal independent family.  
That is, we shall find  $Z\subseteq \omega $  such that, 
whenever  $\Cal X$ and  $\Cal Y$  are 
 disjoint finite subfamilies of  $\Cal I$, then both  
$W(\Cal X,\Cal Y)\cap Z$  and  $W(\Cal X,\Cal Y) - Z$
are infinite, so $\Cal I\cup\{Z\}$ is independent and properly 
includes  $\Cal I$.

Partition  $\Cal I$  as  $\Cal D\cup \Cal E$, where  
$\Cal D = \{D_n\mid n\in \omega \}$  is countably infinite 
 and  $\Cal E$  has (like  $\Cal I$) cardinality smaller than $\dd$.
Write $D_n^0$ and $D_n^1$   for 
$D_n$   and  $\omega -D_n$, respectively.  
For each  $x\in{}^\omega2$, we apply Lemma 22 with
$$
C_n =\bigcap_{k<n}D_k^{x(k)} 
$$
 and
$$
\Cal A = \{W(\Cal X,\Cal Y)\mid\Cal X,\Cal Y 
\text{ finite disjoint subfamilies of  }\Cal E\}. 
$$
 The hypothesis of the lemma is satisfied because  $\Cal I$  is 
independent.  So we obtain  $B_x$   such that
\newline(1)\quad $B_x\subseteq^*\bigcap_{k<n}D_k^{x(k)}$ 
for all $n$, and
\newline(2)\quad $B_x\cap W(\Cal X,\Cal Y)$ is infinite for every  
$W(\Cal X,\Cal Y)\in\Cal A$.
\newline
 Notice that, by (1),
\newline(3)\quad $B_x\cap B_y$ is finite when  $x\neq y$.

Fix two disjoint, countable, dense (in the usual product topology) 
subsets $Q$ and $Q'$ of ${}^\omega2$.  We can remove finitely many
elements from $B_x$, for $x\in Q\cup Q'$, so that
\newline(3*)\quad $B_x\cap B_y=\emptyset$ for $x\neq y$  in  
$Q\cup Q'$.\newline
 (To see this, we use the countability of $Q\cup Q'$ to list the
relevant $B$'s in an $\omega$-sequence, and remove from each one
its (finite, by (3)) intersections with its (finitely many) 
predecessors in the list.)  Notice that (1) and (2) remain true.

Now set
$$
Z =\bigcup_{x\in Q}B_x\quad\text{and}\quad 
Z' =\bigcup_{x\in Q'}B_x, 
$$
 so  $Z$  and  $Z'$  are disjoint, by (3*).  
We shall show that, for any finite  disjoint  
$\Cal X,\Cal Y\subseteq \Cal I$, the intersection  
$W(\Cal X,\Cal Y)\cap Z$ is infinite.  The same reasoning with  
$Q'$  in place of $Q$  will yield that  $W(\Cal X,\Cal Y)\cap Z'$  is 
 infinite, and therefore so is  $W(\Cal X,\Cal Y) - Z$, which will 
complete the proof.

Let finite disjoint $\Cal X,\Cal Y\subseteq\Cal I$  be given.  
As $Q$ is dense in ${}^\omega2$, it 
 contains an $x$ such that, if $D_k\in\Cal X$  (resp. $\Cal Y$), 
then  $x(k)=0$ (resp. 1),  so $D_k^{x(k)}=D$  (resp. $\omega-D$).  
Thus,
$$\align
W(\Cal X,\Cal Y) &= W(\Cal X\cap\Cal E,\Cal Y\cap\Cal E)\cap 
W(\Cal X\cap\Cal D,\Cal Y\cap\Cal D) \\
&= W(\Cal X\cap \Cal E, \Cal Y\cap \Cal E)\cap
\bigcap_{k:D_k\in\Cal X\cup\Cal Y}D_k^{x(k)}\\
&\supseteq W(\Cal X\cap\Cal E,\Cal Y\cap\Cal E)\cap
\bigcap_{k<n}D_k^{x(k)},\\
&\qquad\qquad\text{for sufficiently large }n,\\
&\supseteq^* W(\Cal X\cap\Cal E,\Cal Y\cap\Cal E)\cap B_x,\\
&\qquad\qquad\text{by (1)}.
\endalign
$$
 As  $W(\Cal X\cap\Cal E,\Cal Y\cap\Cal E)\in\Cal A$, 
we have by (2) that its intersection with $B_x$  is infinite.  
We have just seen that this infinite set is almost included in 
 $W(\Cal X,\Cal Y)$, and it is also included in $Z$ because $x\in Q$  
implies $B_x\subseteq Z$. 
 So  $W(\Cal X,\Cal Y)\cap Z$  is infinite, as required. 
\qed\enddemo

\enddocument